\documentclass{scrartcl}

\usepackage[a4paper,margin=2.5cm]{geometry}

\usepackage[utf8]{inputenc}
\usepackage[T1]{fontenc}

\usepackage{graphicx}

\usepackage{amsmath,amssymb,amsthm}
\usepackage{esint} 
\usepackage{enumitem}

\usepackage[dvipsnames,svgnames]{xcolor} 
\usepackage{comment}
\usepackage{caption}
\usepackage{subcaption}

\usepackage{cancel} 

\usepackage[style=alphabetic,maxalphanames=4,maxnames=4]{biblatex}
\AtEveryBibitem{\clearfield{issn}}

\newtheorem{thm}{Theorem}
\newtheorem{prop}[thm]{Proposition}
\newtheorem{lem}[thm]{Lemma}
\newtheorem{cor}[thm]{Corollary}
\theoremstyle{definition}
\newtheorem{Def}[thm]{Definition}
\newtheorem{Def-Thm}[thm]{Definition-Theorem}

\theoremstyle{remark}
\newtheorem{rem}[thm]{Remark}

\newcommand{\dpart}[2]{\frac{\partial #1}{\partial #2}}
\newcommand{\N}{\mathbb{N}}

\newcommand{\Q}{\mathbb{Q}}

\newcommand{\R}{\mathcal{R}}
\newcommand{\C}{\mathbb{C}}

\usepackage[
    colorlinks=true,
    citecolor=NavyBlue,
    linkcolor=NavyBlue,   
    urlcolor=blue
]{hyperref}

\title{Enumeration of plane hypermaps with a mixed boundary II}
\author{Thomas
 Lejeune
  \thanks{Sorbonne Université and Université Paris Cité, CNRS,
    IMJ-PRG, F-75005 Paris, France}}
\date{\today}

\begin{document}
\maketitle

\vspace{-0.9cm}
\begin{abstract}

A plane hypermap with a boundary is a connected graph embedded in the plane with a distinguished face, and whose inner faces are properly face-bicolored. We say that such hypermaps are $k$-alternating if the colors 
of the faces adjacent to the distinguished face alternate at most $2k$ times when going counterclockwise. 
Eynard and Orantin, in \cite{EO2005}, have found a universal expression for the generating functions of $k$-alternating hypermaps as a surprisingly simple function of the generating function of $1$-alternating hypermaps.\\
In this article, we will provide combinatorial formulas for $k$-alternating hypermaps using specific subclasses of hypermaps in which each vertex can either be connected to a previously marked vertex, called \emph{accessible}, or connected to all other vertices in what we call a \emph{blob}.
In particular, this method will allow us to derive the results from \cite{EO2005} from a combinatorial perspective, and we anticipate that it will lead to the computation of asymptotics for hypermaps in future work.\\
This article is naturally a sequel to, and most importantly a generalization of, the results in \cite{BEL26}, which deals essentially with the case of $1$- and $2$-alternating hypermaps.
\end{abstract}

\vspace{-0.5cm}
\tableofcontents

\paragraph{Acknowledgments} I would like to thank Jérémie Bouttier and Bertrand Eynard, who closely supervised the development of this article. Thanks also to Marie Albenque and Gilles Schaeffer for having helped me get to grips with this problem, due to their valuable perspectives on these questions.
\\
This article was not written using any help from an artificial intelligence.

\section{Introduction}

\subsection{General motivations}

Combinatorial maps, that is, graphs embedded in surfaces up to
homeomorphism, have long served as a meeting point for enumerative
combinatorics, statistical mechanics, and random geometry. Since
Tutte's pioneering enumeration of planar
maps~\cite{Tutte_1963}, the subject has grown into a rich theory \cite{Schaeffer2015,Curien2023} with
connections ranging from the KP hierarchy and matrix models in
mathematical physics~\cite{Eynard2016} to scaling limits of random
surfaces~\cite{LeGall2013,Miermont2013}. The present paper contributes
to this circle of ideas by studying \emph{hypermaps}, a natural
generalization of maps obtained by allowing two types of faces, and
by developing bijective techniques adapted to their enumeration.
\\

Hypermaps arise as face-bicolored maps, and their enumeration
according to the degree sequences of both colors is governed by a
two-matrix integral, closely related to the Ising models
on random surfaces~\cite{Kazakov1986,Boulatov1987,tokka2025isingmodelexternalmagnetic,Duits25}. These statistical
mechanics models undergo a continuous phase transition at a specific
critical point, where the underlying random geometry seems to
converge \cite{OAI:Canonical-conformal-limits-of-subcritical-FK-planar-maps-September-24-2026}, in an appropriate scaling limit, to the
Liouville quantum gravity metric with a central charge $c \in (0,1)$ \cite{Ding2023,BP25}.
It then remains to make this abstract convergence more explicit in order to 
make further LQG computations,
and a finer combinatorial understanding of hypermaps, in particular
through explicit bijective constructions, is a natural step toward
that goal.
\\
Beyond their probabilistic and physical motivations, hypermaps also
fit within the broader algebraic framework of free probability, where
the enumeration of maps by face degrees can be recast in terms of free
cumulants and free convolution~\cite{MingoSpeicher2017}; extending
this dictionary to the bicolored setting of hypermaps is an appealing
direction that we leave for future investigation.
\\

In this context, some twenty years ago, Eynard and Orantin \cite{EO2005,Eynard2016} proved a number of elegant formulas
linking hypermaps with a mixed boundary embedded in the plane (known as plane hypermaps) to plane hypermaps with simpler boundary, known as Dobrushin hypermaps.
These formulas, originally derived from a Tutte decomposition as well as a strong algebraic intuition, nevertheless lack
a combinatorial interpretation that would allow us to look even beyond these results.
\\
In the previous article \cite{BEL26}, we studied in greater depth the case
of plane hypermaps with a Dobrushin boundary condition and began to investigate the
case of slightly more complex mixed boundaries. The aim of this article is
to build on the previous work by generalising the observations made in the first. To do so,
we will use the slice decomposition approach \cite{AB2025,BOUTTIER_2014} to split a hypermap with a boundary
containing a marked vertex into a so-called accessible component, containing all vertices in the vicinity of the marked vertex,
as well as other hypermaps with smaller boundaries. We shall therefore focus
on this notion of an accessible component and demonstrate an
algebraic result concerning so-called strongly connected hypermaps, also informally known as blobs,
where all vertices are said to be accessible. To clarify all these concepts,
it is time to define our terminology more rigorously.

\subsection{Lexicon and definitions}

We briefly recall the notations needed to understand the results
that follow, and refer the reader to \cite{Schaeffer2015, BEL26} for a fuller
account of the definitions and preliminary properties underlying
this article.
\\

A planar map is a cellular embedding of a connected graph on the sphere, considered up to deformations,
and a plane map is a planar map with a marked face, which allows it to be
projected stereographically onto the plane so that the marked face becomes
the outer face. This definition allows us to define the faces of the map as well as its corners, which 
we will specifically deal with in Definition \ref{def:corner}. A plane hypermap is then a properly bicolored plane map, that is,
one in which each inner face is assigned a color, black or white, in such a way
that each white face is only adjacent to black faces (or to the outer face) and vice versa. We shall
simply call these plane maps hypermaps in what follows, in order to avoid
redundancy. Note that the outer face has no predefined color. \\

Such a graph then admits a canonical orientation by considering each face as
a cycle, in such a way that white faces are those oriented clockwise and black
faces are those oriented counterclockwise. Only one type of edge remains
unoriented at this point, namely edges that are not incident to any inner face, and so are only incident to the outer face, which
we call bridges. Note that, since we are in the planar case and there is only one distinguished face, a bridge is exactly an edge whose removal disconnects the hypermap, hence its denomination. Bridges will then play a particular role in what follows. For now, they can carry any orientation one likes, however it will 
soon be fixed when we define $k$-alternating hypermap by setting what we shall call the
boundary condition.
\\

\begin{figure}[h!]
    \centering
    \includegraphics[width=0.7\textwidth]{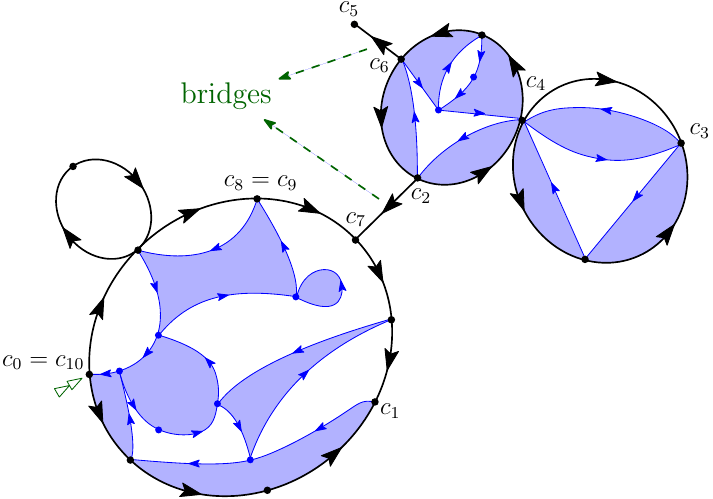}
    \caption{A representation of a plane hypermap $\mathfrak{m}$, where black faces are represented in blue and the boundary is made of black edges. Note the presence of two bridges.
    \\
    By going counterclockwise and counting the marked corners, this is hence the sketch of a $5$-alternating hypermap with weight contribution $w(\mathfrak{m})=\displaystyle\frac{t^{24}\,(t_2^{\circ})^3\,(t_3^{\circ})^3\,(t_4^{\circ})^{2}\,t_5^{\circ}\,t_7^{\circ}\,t_1^{\bullet}\,(t_2^{\bullet})^5\,(t_3^{\bullet})^{3}\,(t_4^{\bullet})^4}{x_1^4\,y_1^4\,x_2^4\,y_2^2\,x_3^4\,y_3^2\,x_4^3\,y_4^2\,x_5^1\,y_5^5}$ to $H_5$.}
    \label{fig:hypermap_generale}
\end{figure}

Figure \ref{fig:hypermap_generale} is a sketch of such a hypermap, where black faces are
represented in blue for aesthetic purposes and the boundary edges are
represented in black.

The set of edges incident to the outer face, in particular containing all
bridges, is called the \emph{boundary} of the hypermap. Along this boundary,
we denote by $[c,c']$ the set of edges and vertices of the boundary appearing,
going counterclockwise, between corners $c$ and $c'$. This is what we call a
boundary interval. A boundary interval is said to be \emph{monochromatic} if all the edges it contains
are oriented in the same direction. Note that this definition in particular
takes bridges into account and that, if the boundary interval contains no bridge, then all faces incident 
to those edges are of the same color.
\\

We now consider a strictly positive integer $k$. A hypermap is said to be
$k$-alternating if its outer face is adjacent to $2k$ marked corners
$c_0,c_1,\dots,c_{2k}=c_0$, appearing counterclockwise, such that the boundary intervals $[c_i,c_{i+1}]$, $i\in[\![0,2k-1]\!]$,
are all monochromatic and oriented alternately counterclockwise (say \emph{black}) then clockwise (say \emph{white}). The choice of colors 
made for these boundary intervals characterizes, as explained above, the colors of the faces every edge of the interval is incident to (except bridges). Note also that these intervals can a priori 
be empty as two consecutive corners could merge into one. We shall say that a boundary is monochromatic if the whole boundary interval $[c_0,c_{2k}]$ is.
\\

\begin{figure}[h!]
    \centering
    \begin{subfigure}[b]{0.45\textwidth}
        \includegraphics[width=\textwidth]{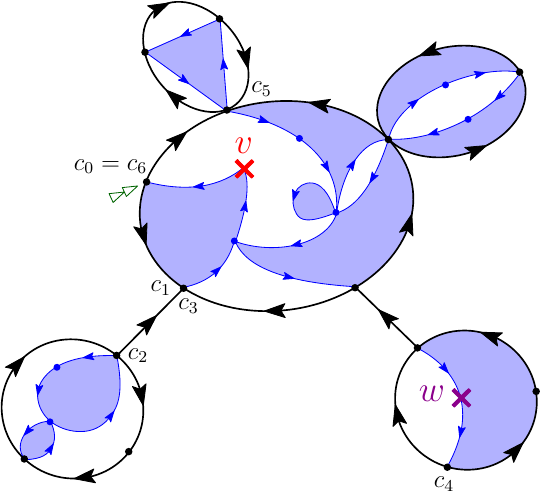}
        \caption{Two possible $3$-alternating pointed hypermaps, accessible if $v$ is marked and not if $w$ is (because of the presence of the bridge near $w$).}
        \label{fig:hypermap_accessible}
    \end{subfigure}
    \hfill
    \begin{subfigure}[b]{0.45\textwidth}
        \includegraphics[width=\textwidth]{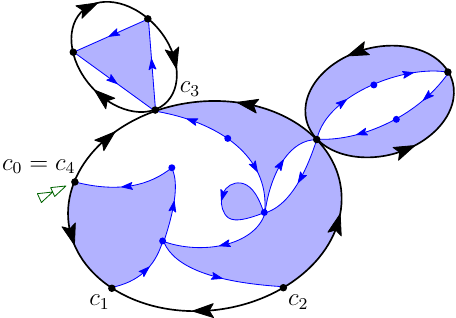}
        \caption{A strongly connected $2$-alternating hypermap. Note that it has no bridge and that any vertex is accessible from any other vertex of the hypermap.}
        \label{fig:hypermap_stronglyconnected}
      \end{subfigure}
    \caption{Examples of specific $k$-alternating hypermaps.}
    \label{fig:hypermap}
\end{figure}

Finally, the main method we shall use in what follows to enumerate
$k$-alternating hypermaps is to consider pointed $k$-alternating hypermaps, that
is, those possessing a marked vertex. Within this class of hypermaps, we shall
mainly consider the class of accessibly pointed $k$-alternating hypermaps. Let us
consider a vertex $v$ of a hypermap: we say $v$ is \textbf{accessible} if every vertex $u$
of the hypermap is connected to $v$ by an oriented path from $u$ to $v$. We then say that a hypermap is \textbf{accessibly pointed} if one of its vertices is marked and the latter is accessible. It is
precisely this notion that will be central throughout this article. We give an example of an accessibly pointed hypermap on Figure \ref{fig:hypermap_accessible}.
\\
Moreover, a $k$-alternating hypermap is said to be \emph{strongly connected} if all of its vertices is accessible. Informally, we
will also call these hypermaps \emph{blobs} as it appears in \cite{AB2025}. See Figure \ref{fig:hypermap_stronglyconnected} for an example of a blob.

\subsection{Weights and generating series}

The method we have chosen to enumerate $k$-alternating hypermaps is that of
formal generating series. We thus consider formal variables $t, t_i^{\circ},
t_i^{\bullet}, i\in\N^*$ such that $t$ is the weight of a vertex, $t_i^{\circ}$
the weight of a white face of degree $i$, and $t_i^{\bullet}$ the weight of a
black face of degree $i$.
\\

As stated earlier, we shall consider boundary intervals of various lengths.
We shall thus mark with a weight $x_i^{-1}$ (resp. $y_i^{-1}$) each edge
appearing on the so-called black boundary interval $[c_{2i-2},c_{2i-1}]$
(resp. the white boundary interval $[c_{2i-1}, c_{2i}]$) for $i\in[\![1,k]\!]$,
using the notations from the previous section. We shall also add a weight
$\prod_{i=1}^kx_i^{-1}y_i^{-1}$, whose purpose is to simplify computations as
well as combinatorial interpretations relative to bridges, as we shall see
later on.\\

To summarize, if we denote by $\mathfrak{h}_k$ the set of $k$-alternating
hypermaps, then the weight of a hypermap $\mathfrak{m}\in\mathfrak{h}_k$ is
written as:
\begin{equation}\label{eq:weight}
w(\mathfrak{m}):=t^{V(\mathfrak{m})}\prod_{i\in\N^*}(t_i^{\circ})^{F_i^{\circ}(\mathfrak{m})}\prod_{j\in\N^*}(t_j^{\bullet})^{F_j^{\bullet}(\mathfrak{m})}\prod_{i=1}^kx_i^{-\ell_i^{\bullet}-1}\prod_{j=1}^ky_j^{-\ell_j^{\circ}-1},
\end{equation}
where $V(\mathfrak{m}), F_d^{\circ}(\mathfrak{m})$ and $F_d^{\bullet}(\mathfrak{m})$ are respectively the number of vertices, white faces of degree $d$ and black faces of degree $d$ of $\mathfrak{m}$, for $d\in\N^*$. Also, for $i\in[\![1,k]\!]$, $\ell_i^{\circ}(\mathfrak{m})$ and $\ell_i^{\bullet}(\mathfrak{m})$ are respectively the number of edges appearing on the boundary interval $[c_{2i-2},c_{2i-1}]$ and $[c_{2i-1},c_{2i}]$ of $\mathfrak{m}$. As explained before, these 
lengths can equal zero if the boundary interval is empty.
\\

We then denote by $H_k$ the generating function of $k$-alternating hypermaps:
\begin{equation}
  \label{eq:Hkdef}
  H_k(x_1,\ldots,x_k;y_1,\ldots,y_k; t, (t_i^{\circ})_{i\geqslant1},(t_j^{\bullet})_{j\geqslant1}) := \delta_{k,1} + \sum_{\mathfrak{m}\in\mathfrak{h}_k} w(\mathfrak{m}),
\end{equation}
where one will notice the appearance of an additional term $\delta_{k,1}$,
which will allow us to distinguish the presence, or absence, of a
$1$-alternating hypermap in our future reasoning. For the sake of clarity, we
shall omit writing the weights of faces or vertices as arguments, and we shall
even favor the notation $H(1\to k)$ instead of
$H_k(x_1,\ldots,x_k;y_1,\ldots,y_k)$ whenever it proves easier to manipulate. See Figure \ref{fig:hypermap_generale} for an example of a $k$-alternating hypermap and its associated weight.
\\

\begin{figure}[h!]
    \centering
    \begin{subfigure}[b]{0.45\textwidth}
        \includegraphics[width=\textwidth]{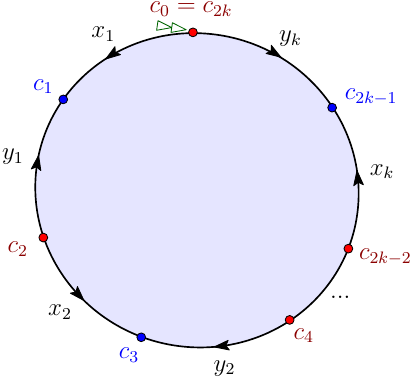}
        \caption{Generic sketch of a $k$-alternating hypermap counted by $H_k$. Inner faces are blurred and marked corners are represented.}
        \label{fig:schema_sans_rien}
    \end{subfigure}
    \vfill
    \begin{subfigure}[b]{0.45\textwidth}
        \includegraphics[width=\textwidth]{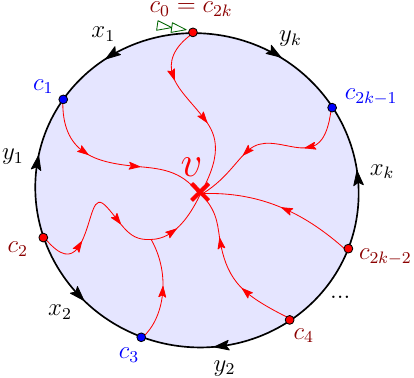}
        \caption{Generic sketch of an accessibly pointed $k$-alternating hypermap counted by $A_k$. The marked vertex is denoted as $v$ and only the oriented paths from all marked corners are represented on the figure, see Lemma \ref{lem:accessible} for an explanation of this choice.}
        \label{fig:schema_sans_rien_accessible}
    \end{subfigure}
    \hfill
    \begin{subfigure}[b]{0.45\textwidth}
        \includegraphics[width=\textwidth]{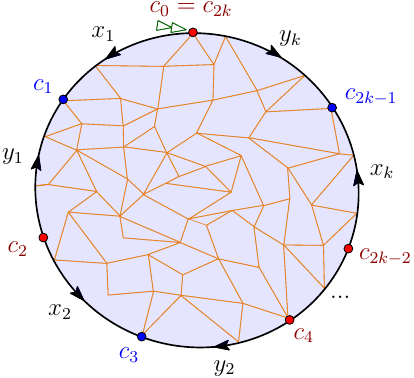}
        \caption{Generic sketch of a strongly connected $k$-alternating hypermap counted by $B_k$. We represented the inner connections with light orange here, inspired by the slime mold \emph{Physarum polycephalum} (also called "blob").}
        \label{fig:schema_sans_rien_stronglyconnected}  
    \end{subfigure}
    \caption{Sketches of specific families of $k$-alternating hypermaps}
    \label{fig:schema_sans_rien2}
\end{figure}

We finally recall the introduction of the sets of accessibly pointed and
strongly connected $k$-alternating hypermaps defined above. The generating
functions associated with these sets will be denoted respectively by $A_k$ and
$B_k$ in what follows and defined in the same way as $H_k$, and will not carry the additional term $\delta_{k,1}$.
In our combinatorial decomposition, we will use generic sketches to represent hypermaps contributing to $H_k, A_k$ and $B_k$, see Figure \ref{fig:schema_sans_rien2}.
These series are the ones that will lie at the heart of our
future reasoning and theorems, as we shall now go on to explain.

\subsection{What has been proved in the first article...}

\subsubsection{...about the methodology}

Let us recall that this article is the sequel to \cite{BEL26}, where we
introduced the first ideas that led to the generalizations we shall present
here. The previous article focused on creating a new way of enumerating
hypermaps using the notion of accessibility. We had then re-derived various
formulas appearing in \cite{AB2025} and \cite{Eynard2016} concerning hypermaps
with a monochromatic boundary or $1$- and $2$-alternating boundaries, in
addition to finding new ones. We therefore refer the reader to these
precursor articles, which provide the necessary understanding of what led to
this work.
\newpage

The main observation from the previous article, which will guide the structure of
this one, is the following: there exists a specific family of
hypermaps, the one of accessibly pointed hypermaps, which allows any pointed hypermap to be
decomposed according to its \emph{accessible component}, that is the set of
vertices connected to the marked vertex, together with other hypermaps with a
smaller number of marked corners on the boundary. Moreover, the algebraicity of Tutte's method for planar maps
ensures that we can integrate the differential equation thus obtained by distinguishing a vertex, equivalent to the operator $\dpart{}{t}$ for generating series, in
order to find algebraic formulas for these hypermaps. In this article, we
shall substantiate this thesis for the case of general $k$-alternating
hypermaps, so no longer just for $k=1$ or $2$ as in \cite{BEL26} but for all $k\geqslant1$.

\subsubsection{...about the results}\label{sec:resultBEL}

Let us consider the following formal Laurent power series
\begin{equation}\label{eq:defxynew}
  x(z):=\sum_{i\geq -1} a_i z^{-i}\quad\text{ and }\quad y(z):=\sum_{i \geq -1} b_i z^i
\end{equation}
where the $a_i$'s and $b_i$'s are elements of $\R$, which is defined as $\Q[\![t,t_1^\circ,t_2^\circ,\ldots,t_1^\bullet,t_2^\bullet,\ldots]\!]$, the set of formal power series 
in the physics variables $t,t_1^\circ,\ldots,t_1^\bullet,\ldots$ with rational coefficients. They can be
determined recursively by the conditions
\begin{equation}
  \label{eq:aibieq}
  \begin{split}
    a_i = t \delta_{i,-1} + \sum_{d \geq 1} t^\bullet_d [z^{-i}] y(z)^{d-1} \qquad (i \geq -1)\\
    b_{-1} = 1, \qquad  b_i = \sum_{d \geq 1} t^\circ_d [z^i] x(z)^{d-1} \qquad (i \geq 0).
  \end{split}
\end{equation}

These Laurent power series can be interpreted both as generating series of hypermap slices,
studied in \cite{AB2025}, and, in the case of bounded face degrees, as a
rational parametrization of the natural spectral curve associated with our model, see \cite{Eynard2003,Eynard2016}. We refer the
interested reader to these articles; the other reader may simply note that these quantities are well understood. We won't need to understand deeply their interpretation in our future reasoning.
\\

In fact, we have shown in \cite{BEL26} that the
generating function of accessibly pointed $k$-alternating hypermaps can be expressed thanks to $x(z)$ and $y(z)$:
\begin{equation}\label{eq:ak}
\begin{split}
    A_k(x_1,\ldots,x_k;y_1,\ldots,y_k) &= \sum_{\ell_1^\bullet,\ell_1^\circ,\ldots,\ell_k^\bullet,\ell_k^\circ \geq 0} \frac{[z^0] x(z)^{\ell_1^\bullet+\cdots+\ell_k^\bullet} y(z)^{\ell_1^\circ+\cdots+\ell_k^\circ}}{x_1^{\ell_1^\bullet+1}y_1^{\ell_1^\circ+1}\cdots x_k^{\ell_k^\bullet+1}y_k^{\ell_k^\circ+1}}\\
    &=[z^0] \frac{1}{\prod_{i=1}^k (x_i-x(z))(y_i-y(z))},
\end{split}
\end{equation}
so we can claim we know $A_k$ explicitly for every $k\in\N^*$.
\\

We have then been able to determine in particular two relations between $H_1, H_2$ in terms
of $A_1,A_2$, namely:

\begin{equation}\label{eq:dH1/dt}
\begin{split}
    \dpart{H_1(x_1;y_1)}{t}=H_1(x_1;y_1)\, A_1(x_1;y_1),
\end{split}
\end{equation}
and
\begin{equation}\label{eq:dH2/dt}
\begin{split}
\dpart{H_2(x_1,x_2;y_1,y_2)}{t}
&= A_2(x_1,x_2;y_1,y_2)\,H_1(x_1;y_1)\,H_1(x_2;y_2)\\
&+ \left(A_1(x_1;y_2)+A_1(x_2;y_1)\right)\,H_2(x_1,x_2;y_1,y_2).
\end{split}
\end{equation}
We had thus integrated the latter Equation \eqref{eq:dH2/dt} using \eqref{eq:dH1/dt} and \eqref{eq:ak} in order to
find the elegant relation, appearing in \cite{Eynard2003}:
 \begin{equation}\label{eq:h2}
    H_2(x_1,x_2;y_1,y_2) = \frac{H_1(x_1;y_1)H_1(x_2;y_2)-H_1(x_1;y_2)H_1(x_2;y_1)}{
      (x_1-x_2)(y_1-y_2)}.
  \end{equation}
We shall see in this article a generalization to any $k$ of Equations \eqref{eq:dH1/dt}
and \eqref{eq:dH2/dt} as our Theorem \ref{th:dHk/dt}, whose integration will allow us to find a more general
form of \eqref{eq:h2}, stated in Corollary \ref{cor:hk}.
\\

Regarding the series $B_2$, we had also proved the relation:
 \begin{equation}\label{eq:b2}
    B_2(x_1,x_2;y_1,y_2) = - \frac{H_1(x_1;y_1)^{-1}H_1(x_2;y_2)^{-1}-H_1(x_1;y_2)^{-1}H_1(x_2;y_1)^{-1}}{
      (x_1-x_2)(y_1-y_2)},
  \end{equation}
  which we shall generalize here in Theorem \ref{th:bk}. We also recall the relation for $k=1$:
 \begin{equation}\label{eq:h1b1}
    H_1(x_1;y_1) = \frac{1}{1-B_1(x_1;y_1)},
  \end{equation}
  proved in \cite{AB2025} and \cite{BEL26} and which will be useful to us in Section \ref{sec:bk}.

\subsection{Goals and outline of the article}
As explained previously, the goal of this article is twofold. It is to use
the accessibility method in order to derive new (differential) equations on the series
$H_k$, and then to integrate these relations to recover the theorems of
Chapter 8 of \cite{Eynard2016}, first proved in \cite{EO2005}, thereby demonstrating the strength of this
method. In particular, we shall prove the following theorem:
\begin{thm}\label{th:dHk/dt}
    The decomposition of a pointed $k$-alternating hypermap into accessible components yields the relation
 \begin{equation}\label{eq:dHk/dt}
    \dpart{H(1\to k)}{t}=\sum_{p=1}^k\sum_{1\leqslant i_1<\dots<i_p\leqslant k}A_p(x_{i_1},\dots,x_{i_p};y_{i_1-1},\dots,y_{i_p-1})\prod_{j=1}^{p}H(i_j\to i_{j+1}-1),
  \end{equation}
  where we introduce the shorthand notations $i_{p+1}:=i_1, y_0:=y_k$ and:
  $$H(a\to b):=\left\{\begin{array}{lll}
&H_{b-a+1}(x_a,\dots,x_b;y_a,\dots,y_b)\,\,\,\text{if $a\leqslant b$},&\\
&H_{k+1+b-a}(x_a,\dots,x_k,x_1,\dots,x_b;y_a,\dots,y_k,y_1,\dots,y_b)\,\,\,\text{if $a>b$},
\end{array}\right.$$
which allows us to preserve the rotational invariance of $H_k$ with respect
to its variables $x_i$ and $y_i$ within Formula \eqref{eq:dHk/dt}.
\end{thm}

We shall thus prove this theorem in Section \ref{sec:dHk/dt}. In the
following Section \ref{sec:hk} and \ref{sec:2.4}, we shall deduce from this theorem the
following corollary.
\begin{cor}[\cite{EO2005}]\label{cor:hk}
There exist some coefficients $(C_{\sigma})_{\sigma\in\mathfrak{S}_k}$ indexed by the set of permutations of $[\![1,k]\!]$ which are rational functions of the $x_i$'s and $y_i$'s and which are universal, in the sense that they do not depend on the $t_i^{\circ},t_i^{\bullet}$ nor on
$t$, such that
\begin{equation}\label{eq:hkh1}
H_k(x_1,\dots,x_k;y_1,\dots,y_k)=\sum_{\sigma\in\mathfrak{S}_k}C_{\sigma}(x_1,\dots,x_k;y_1,\dots,y_k) \prod_{i=1}^kH_1(x_i;y_{\sigma(i)}).
\end{equation}
\end{cor}
Furthermore, we shall see that $C_{\sigma}\neq0$ if $\sigma$ is a so-called a plane permutation.
The definition of a plane permutation is strongly related to the concept of non-crossing chord diagram and will be given in Definition \ref{def:planaire}. We will also give 
other proofs of a few ways to compute recursively the $C_{\sigma}$ in Section \ref{sec:hk} and Section \ref{sec:bk}.
\begin{rem}
In fact, these coefficients $C_{\sigma}$ are unique. This is a property that we won't prove explicitly in this article but the method will be the same as for the next theorem, starting from 
Theorem \ref{th:h1param}.
\end{rem}

We shall then generalize Relation \eqref{eq:b2} in Section
\ref{sec:bk}, which will yield the following result.
\begin{thm}\label{th:bk}
Let $k\geqslant2$ be an integer. The generating series $B_k$ of strongly connected $k$-alternating hypermaps
is given by the relation:
$$
B_k(x_1,\dots,x_k;y_1,\dots,y_k)=-\sum_{\sigma\in\mathfrak{S}_k}C_{\sigma}(x_1,\dots,x_k;y_1,\dots,y_k) \prod_{i=1}^kH_1(x_i;y_{\sigma(i)})^{-1}.
$$
where the $C_{\sigma}$ are the same as those defined in Corollary \ref{cor:hk}.
\end{thm}

\begin{rem}
$\bullet$ Corollary \ref{cor:hk} and Theorem \ref{th:bk} are indeed generalizations of formulas \eqref{eq:h2} and \eqref{eq:b2} for $k=2$ found in \cite{BEL26} by taking $$C_{\mathrm{Id}}(x_1,x_2;y_1,y_2):=\frac{1}{(x_1-x_2)(y_1-y_2)}=-C_{(12)}(x_1,x_2;y_1,y_2).$$

$\bullet$ Using Equation \eqref{eq:h1b1}, one can finally find a direct closed relation between $k$-alternating blobs and $1$-alternating blobs. Furthermore, if one defines $\mathcal{B}_k:=\delta_{k,1}-B_k$, then one can write Theorem \ref{th:bk} as 
\begin{equation}
\mathcal{B}_k(x_1,\dots,x_k;y_1,\dots,y_k)=\sum_{\sigma\in\mathfrak{S}_k}C_{\sigma}(x_1,\dots,x_k;y_1,\dots,y_k) \prod_{i=1}^k\mathcal{B}_1(x_i;y_{\sigma(i)}),
\end{equation}
for all $k\geqslant1$, and this is exactly the same structure appearing in Corollary \ref{cor:hk} !
\end{rem}

\section{A formula for $k$-alternating hypermaps}

We now turn to the heart of the matter. In this section we prove
Theorem \ref{th:dHk/dt} together with its Corollary \ref{cor:hk}.

\subsection{A few reminders on accessibility}

In order to prove Theorem \ref{th:dHk/dt} — and thus to see how the
generating functions $A_p$ of accessibly pointed hypermaps naturally arise —
we must first revisit a few properties of the directed graph
constructed on our hypermap.

\begin{Def}\label{def:corner}
Recall that a corner is defined as the angular opening between two edges
appearing consecutively when turning counterclockwise around a vertex.
\\
We define two types of corners appearing on the boundary of an alternating
hypermap.
\\

An \emph{inward corner} is a marked corner, of the boundary of
an alternating hypermap, with two ingoing outer edges, it is forcefully of odd index.
The marked corners with two outgoing outer edges are called \emph{outward corners}, they have even indices. Finally, a marked 
corner which is inward or outward will be called \emph{generic}.\\
Figure \ref{fig:corner} gives a sketch of these kind of corners.

\begin{figure}[h!]
    \centering
    \begin{subfigure}[b]{0.4\textwidth}
        \includegraphics[width=\textwidth]{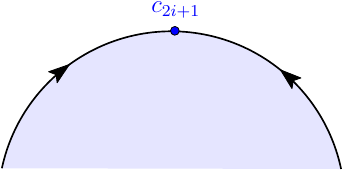}
        \caption{An inward corner, hence with an odd index and represented in blue}
    \end{subfigure}
    \hspace{2cm}
    \begin{subfigure}[b]{0.4\textwidth}
        \includegraphics[width=\textwidth]{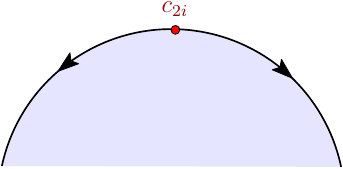}
        \caption{An outward corner, hence with an even index and represented in red}
    \end{subfigure}
    \caption{Sketches of generic corners.}
    \label{fig:corner}
\end{figure}
\end{Def}

\begin{rem}\label{rem:inward_outward}
As said in the definition, inward corners have odd indices
whereas outward corners have even indices. The only
counterexamples preventing from making this a general definition
represent the cases where multiple marked corners are at the same
location. Figure \ref{fig:hypermap_generale} gives an example of two marked corners, $c_8$ and $c_9$, that are respectively not outward and 
not inward.
\end{rem}

This distinction between marked corners will be important in the following reasoning because we really need to consider outward and inward corners to deal with accessibility more than 
other ones. A first example is the following lemma.

\begin{figure}[h!]
    \centering
    \begin{subfigure}[b]{0.4\textwidth}
        \includegraphics[width=\textwidth]{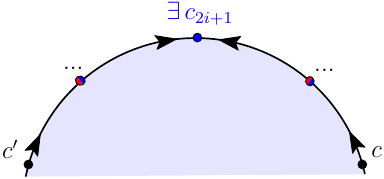}
        \caption{On any inward boundary interval lies an inward corner}
    \end{subfigure}
    \hspace{2cm}
    \begin{subfigure}[b]{0.4\textwidth}
        \includegraphics[width=\textwidth]{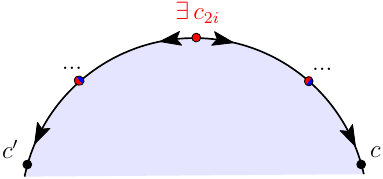}
        \caption{On any outward boundary interval lies an outward corner}
    \end{subfigure}
    \caption{Illustration of Lemma \ref{lem:boundaryinterval}, here bicolored vertices and associated dots stand for a sequence of irrelevant corner lying on the boundary.}
    \label{fig:boundary_interval}
\end{figure}

\begin{lem}\label{lem:boundaryinterval}
Consider a $k$-alternating hypermap and a boundary interval $[c,c']$ of length greater than or equal to $2$ of this hypermap. Then one has the following distinction.

\begin{itemize}
\item If the first edge of $[c,c']$ is oriented counterclockwise and the last one is oriented clockwise, then an inward corner is incident to $[c,c']$. We call it an \emph{inward boundary interval}.
\item If the first edge of $[c,c']$ is oriented clockwise and the last one is oriented counterclockwise, then an outward corner is incident to $[c,c']$. We call it an \emph{outward boundary interval}.
\end{itemize}
\end{lem}
Figure \ref{fig:boundary_interval} pictures this statement that can be understood as an intermediate value theorem.

\begin{proof}
Say $[c,c']$ is an outward boundary interval, the same reasoning apply for the other case. We will use an induction on the length $\ell$ of this interval to conclude. The case $\ell=2$ being trivial, let us assume 
that the property holds for any outward boundary interval of length $\ell-1, \ell\geqslant3$.
\\

Consider the second edge of the interval, if it is oriented counterclockwise then the second corner of $[c,c']$ must be an outward corner. If not, the edge is oriented clockwise and the induction applies on the remaining part of the boundary interval, which concludes.
\end{proof}

Recall that a pointed hypermap is said to be accessibly pointed if its marked vertex is accessible, meaning that it can be
reached from every other vertex of the hypermap by an oriented path. From now on, for simplicity, we
shall identify the properties of the marked corners $c_i$ with those of the vertices
on which each of them lies. In particular, we may speak of a path between a
vertex and a corner. We can now combine these definitions to state the
following lemma.

\begin{lem}\label{lem:accessible}
Let $k$ be a strictly positive integer. A vertex $v$ of a $k$-alternating
hypermap is accessible if, and only if, all the inward corners of the
boundary can reach $v$ by a directed path.
\\

Furthermore, as showed in \cite{BEL26}, a $k$-alternating hypermap has no bridge if, and only if, all
its vertices are accessible. Such a hypermap is, we recall, called strongly
connected.
\end{lem}

\begin{proof}
The direct implication of the first proposition is obvious; let us turn to
the converse. We consider a vertex $v$ that can be reached from every inward
corner of the boundary of the hypermap.
\\

Let us note that every vertex of the boundary that is not an inward corner
can reach one of its neighbors by definition. Thus every vertex of the boundary can reach
an inward corner by iterating this process. It now suffices to show that
every vertex of the hypermap can reach the boundary in order to conclude,
since all inward corners can reach $v$.
\\

We thus consider an arbitrary vertex $w$ of the hypermap, which we connect to
the boundary without taking into account the orientation of the edges. This
is indeed possible since our graph is connected. Let us then consider
consecutively the edges of the path formed from $w$. If an edge is in the
right direction, we keep it; otherwise, we go around an inner face incident
to this edge, which is oriented as a cycle. If it is not incident to an
inner face, it must be incident to the outer face. In other words, it is a
boundary edge, hence $w$ is related to the boundary, which concludes the proof. We refer to Figure
\ref{fig:hypermap_accessible} to visualize the path and we recall Figure \ref{fig:schema_sans_rien_accessible}
where we only emphasized the path from the marked corners to ensure accessibility.
\\

$\hookrightarrow$ We recall here the proof from \cite{BEL26}. \\
First, the existence of a bridge prevents one of the vertices of the bridge from reaching the
other according to its orientation, since erasing a bridge disconnects the hypermap.
Indeed, if one can find a path from the first to the latter, then adding the bridge to the end 
of the path creates a cycle $\mathcal{C}$. But $\mathcal{C}$ contains a bridge, which is only incident to the outer face. And the outer face is a connected component of the plane minus our graph, and so 
cannot be separated by $\mathcal{C}$.\\
Thus, in particular, we cannot connect all vertices of a hypermap containing a bridge to any other ones. \\

Conversely, let us connect two vertices of the hypermap as before, that is,
without taking into account the orientations. We then consider the edges
along the path that are not in the right direction and that are not incident
to any inner face, since these are the ones that could obstruct the creation
of a directed path, as explained earlier. Such an edge can thus only be
incident to the outer face, and is therefore a bridge. However, there can be
no bridge by hypothesis, which concludes the proof.

\end{proof}

The strength of this lemma is that it reduces the notion of accessibility,
which a priori seemed rather restrictive, to simply checking whether each
marked boundary (inward) corner can reach a previously marked vertex, or not. Also, one
intuitively notices that it is therefore the outward corners that are the most likely to being able to reach any vertex. This intuition is detailed in the following lemma.

\begin{lem}\label{lem:au_moins_un}
For a pointed $k$-alternating hypermap, either at least one of its \underline{outward} corners
can reach the marked vertex, or there is no outward corner at all, hence the boundary is monochromatic.
\end{lem}

\begin{proof}
Figure \ref{fig:au_moins_un} will give a sketch of this proof.
\newpage 

Assume the hypermap $\mathfrak{m}$ is not monochromatic and that any outward corner cannot reach the marked vertex. Let $c_0,c_1,\dots,c_{2k-1},c_{2k}=c_0$
denote the sequence of marked corners of the hypermap on its boundary, and
let $v$ denote the marked vertex. \\

Let us prove first that there exists an outward corner on the boundary of $\mathfrak{m}$. Indeed, the boundary is not monochromatic, so there exist two edges on the boundary that are oriented differently. They are then the starting and ending edges of an outward boundary 
interval, which leads to the existence of an outward corner, say it is $c_0$. \\

Thus $c_0$, which is an outward corner, cannot reach $v$. Then this means that there exists a bridge $b_1$ oriented from the marked
vertex towards $c_0$ on all paths from $c_0$ to $v$; let us take the
last of these. This bridge thus separates our hypermap into two parts (excluding $b_1$ itself), one containing $c_0$ and the other containing $v$, let us call the latter $\mathfrak{m}'$. Note that $\mathfrak{m}'$ is a part of $\mathfrak{m}$ since we haven't erased $b_1$.
\\
Now consider the boundary interval of $\mathfrak{m}$ starting from the ending vertex of $b_1$, then going counterclockwise around $\mathfrak{m}'$, and ending at the ending vertex of $b_1$, call it $[c,c']$. Then $[c,c']$ is obviously an outward boundary interval hence Lemma \ref{lem:boundaryinterval} applies and there 
exists an outward corner lying on this interval, so on the boundary of $\mathfrak{m}$. Let us take the first one appearing when turning counterclockwise and denote it as $c_{2i}$.
\\

Since $c_{2i}$ is an outward corner, it thus
cannot reach $v$, and there therefore still exists a bridge, say $b_2$,
oriented from $v$ towards it. Note that $b_2$ has to be oriented differently than $b_1$ since it cannot be on any path from $c_0$ to $v$ since $b_1$ was chosen to be the last one. Then consider the boundary interval starting from $c$, going counterclockwise on $\mathfrak{m}'$ and ending at $c''$, the first corner incident to the ending point of $b_2$ in this path. By construction, $[c,c'']$ is an outward boundary interval and hence there must exist, by Lemma \ref{lem:boundaryinterval}, an 
outward corner lying on $[c,c'']$. Then it should belong to $[c,c_{2i}]$ and $\mathfrak{m}'$ but it is not $c_{2i}$ as the latter is the first encountered in $[c,c']$. This is then absurd, and hence concludes.

\begin{figure}[h!]
    \centering
    \includegraphics[width=\textwidth]{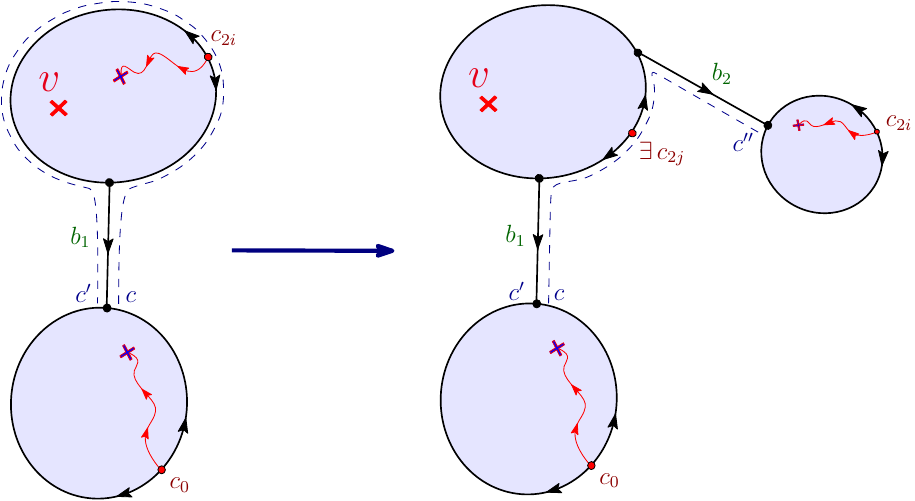}
    \caption{Sketch of the proof of Lemma \ref{lem:au_moins_un}. The different boundary intervals stated are dashed.}
    \label{fig:au_moins_un}
\end{figure}

\end{proof}

\begin{rem}\label{rem:b1_b2}
Note that, in this proof, the fact that $c_0$ and $c_{2i}$ are outward was not used at all. Hence our argument show that if two (not necessarily outward) consecutive marked corners with even indices
cannot reach the marked vertex, then there must exist a common (last) bridge that prevents both of them from reaching it. 
\end{rem}

We are now ready to enumerate pointed $k$-alternating hypermaps and prove Theorem \ref{th:dHk/dt}.

\subsection{Relation between pointed $k$-alternating hypermap and accessibility} \label{sec:dHk/dt}

\begin{figure}
    \centering
    \includegraphics[width=1\textwidth]{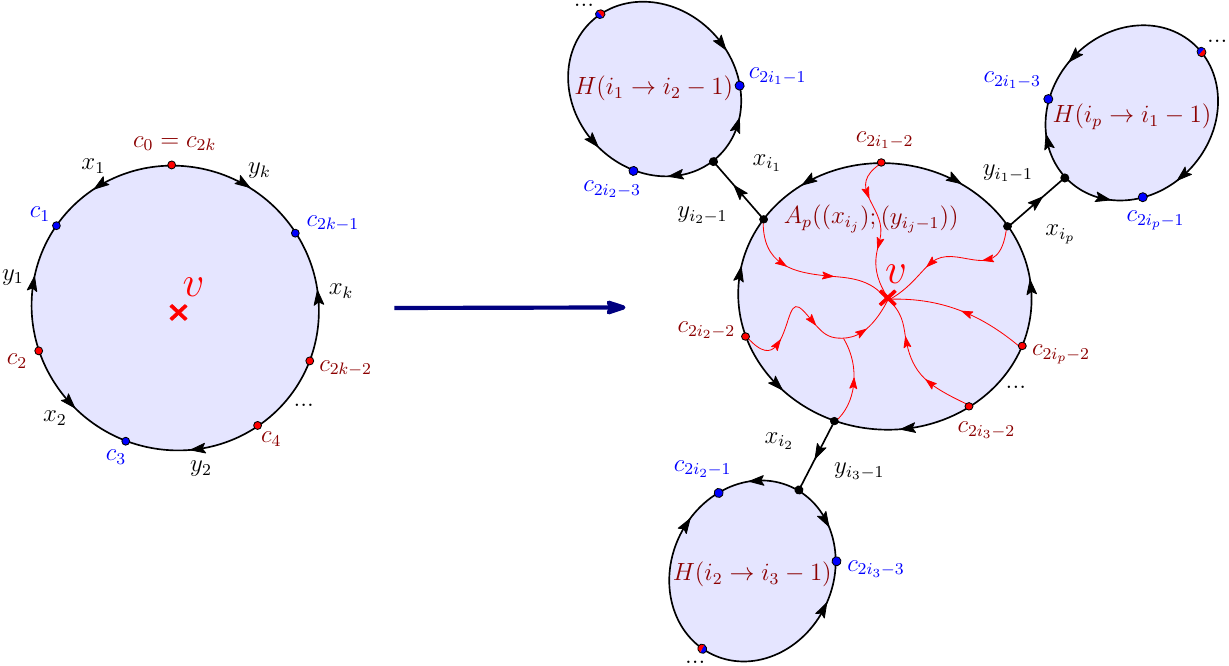}
    \caption{Sketch of the proof of Theorem \ref{th:dHk/dt}. Here a pointed $k$-alternating hypermap 
    is decomposed as a central accessible component containing all the vertices that can reach the mark $v$, and smaller satellite hypermaps 
    attached to its inward corners by a bridge preventing any of its vertices to reach $v$.}
    \label{fig:preuve_totale}
\end{figure}

The whole purpose of this subsection is to prove Theorem \ref{th:dHk/dt}. The main formula given in this 
theorem is a decomposition, which will be illustrated in Figure \ref{fig:preuve_totale}, of a pointed $k$-alternating hypermap (on the left) as a \emph{central} accessibly pointed hypermap which is made 
of all the vertices that can reach the marked vertex, and other \emph{satellite hypermaps} that are attached to it thanks to an outward bridge.
Then let us consider a pointed $k$-alternating hypermap, which does not have a monochromatic boundary (this case will be treated in Remark \ref{rem:monochromatic}), and
where $v$ will be the pointed vertex. Let us further denote, as in the introduction,
$c_0,c_1,\dots,c_{2k-1},c_{2k}=c_0$ the sequence of its marked corners on its
boundary.
\\

Say that there exists $p\in [\![1,k]\!]$ distinct outward corners that can reach $v$, 
denoted as $c_{2i_1-2},\dots,$\\ $c_{2i_p-2}$ where 
$1\leqslant i_1<\dots<i_p\leqslant k$. Note that $p$ cannot equal 
$0$ by Lemma \ref{lem:au_moins_un} and since our hypermap has not a monochromatic boundary.
Then we will show that the accessibly pointed component surrounding $v$ 
contributes by $A_p(x_{i_1},\dots,x_{i_p},y_{i_1-1},\dots,y_{i_p-1})$, 
whose outward corners are exactly $c_{2i_1-2},\dots, c_{2i_p-2}$, and all 
the satellites attached to this components are counted in the product
 $\prod_{j=1}^pH(i_j\to i_{j+1}-1)$, as written in the right-hand side 
 of Equation \eqref{eq:dHk/dt}. This will conclude the proof of Theorem 
 \ref{th:dHk/dt}.
\\
    
To do so, we will first prove that all the corners belonging to a boundary interval strictly 
between two consecutive outward corners that can reach $v$ (let us call them $c_{2i}$ and $c_{2j}$) also all belong to 
a satellite alternating hypermap attached to the boundary interval $[c_{2i},c_{2j}]$ by a common bridge preventing them from reaching $v$.
\\

\begin{figure}
    \centering
    \includegraphics[width=\textwidth]{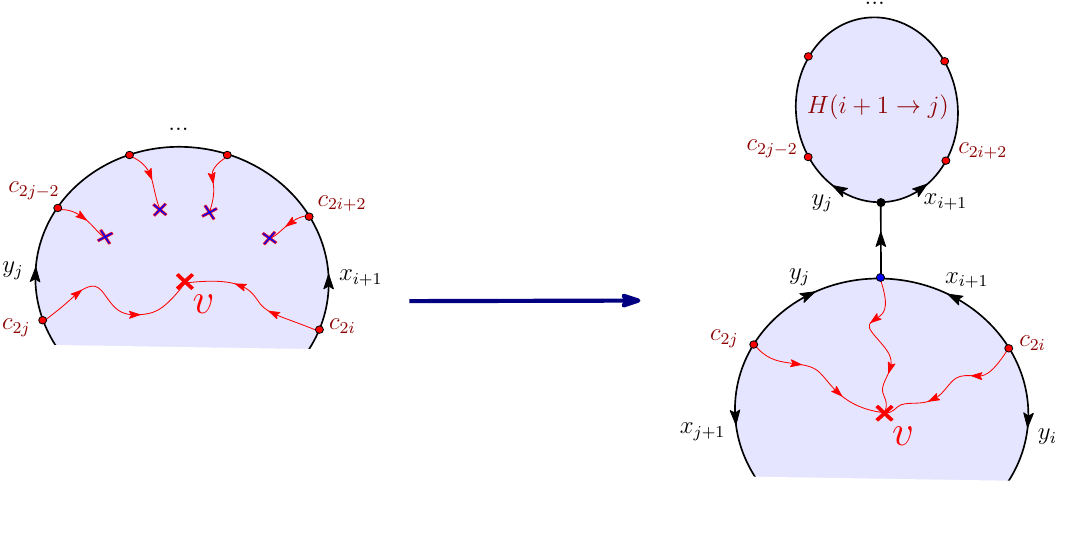}
    \caption{Creation of an isolated component with outward corners that cannot reach $v$. The main information 
    is the presence of a common last bridge that prevents all intermediate outward corners from reaching $v$.}
    \label{fig:c2i_c2j}
\end{figure}

Let us now consider two outward corners that \underline{cannot} reach $v$, if such
corners exist on $[c_{2i},c_{2j}]$ (this is equivalent to assume $j\neq i+1$). Let us also assume them to be consecutive, that is, their
indices are spaced by $2$. As described in the proof of Lemma
\ref{lem:au_moins_un}, the last bridge preventing them from reaching $v$ must
separate the hypermap into two parts: one containing these two
corners and the other containing $v$, as illustrated in Figure \ref{fig:c2i_c2j}. Recall that, if our two corners were
not in the same component, we have shown in the proof of Lemma \ref{lem:au_moins_un} that there necessarily existed an
outward corner between these two corners, and thus they could not be
consecutive. This result can be generalized to any sequence of consecutive corner that cannot 
reach the marked vertex. Indeed, consider the last bridge $b$ that prevents two consecutive outward corners 
(call one of them $c$ and the last of them $c'$ by going counterclockwise on the boundary) to reach $v$, then consider the next outward corner $c''$ on the boundary and assume that
$b$ does not completely prevent $c''$ to reach $v$. Then $c'$ and $c''$ are two consecutive outward 
corners that cannot reach $v$ so, by Lemma \ref{lem:au_moins_un} and Remark \ref{rem:b1_b2}, there must exist a common bridge $b'$ that prevents 
them from reaching $v$. But $b$ should have been the last bridge that prevents $c$ and $c'$ from reaching $v$, which is absurd given the existence of $b'$.
\\
What this intermediate result allows us to say is that all the
outward corners (and even all the even-indexed corners by Remark \ref{rem:b1_b2}) appearing strictly between two consecutive outward corners that can
reach $v$ belong to the same component, since the last bridge separating them from $v$ is the same for
all of them. By detaching this component from the overall hypermap, we thus
obtain a pointed alternating hypermap, the one containing $v$, and a general
$k$-alternating hypermap of weight $H(i+1\to j)$. Note that this detachment is bijective since we only need to mark 
the new hypermap at the starting vertex of the erased bridge and the indices taken into account in $H(i+1\to j)$ makes it easy to know where it has to be attached (at the inward corner created between $c_{2i}$ and $c_{2j}$ by erasing the bridge).
\\

\begin{figure}
    \centering
    \includegraphics[width=\textwidth]{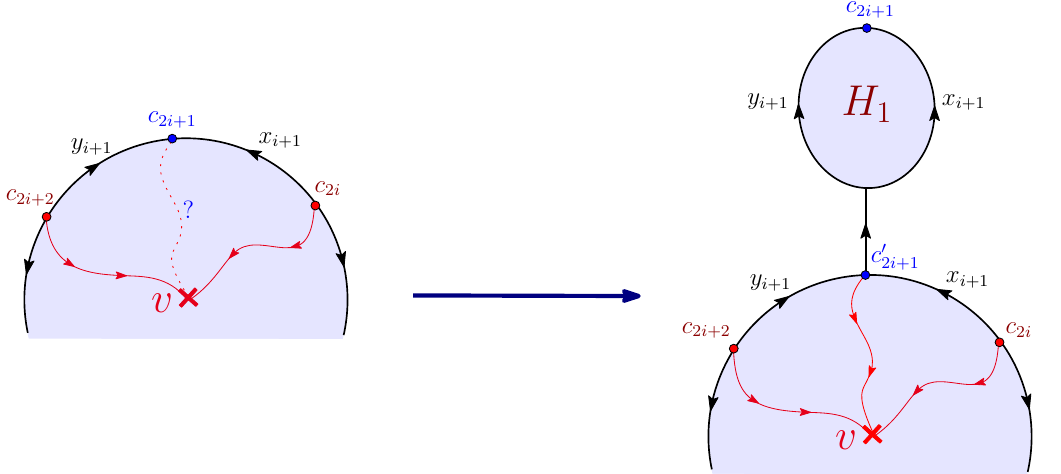}
    \caption{Case of two consecutive outward corners related to $v$. The inward corner $c_{2i+1}$ can or cannot reach $v$, 
    which leads to the existence, or not, of a bridge on its path toward $v$, hence it belongs to a $1$-alternating hypermap. In the case where $c_{2i+1}$ can reach $v$, we set $c_{2i+1}'=c_{2i+1}$ and there is no extra $1$-alternating hypermap as the one in the sketch.}
    \label{fig:coins_adjacents}
\end{figure}

We can iterate this process by going around all the couples of consecutive outward corners $c_{2i},c_{2j}$ that can 
reach $v$. But one case remains to be treated: what is there is no even-indexed corner between $c_{2i}$ and $c_{2j}$ that cannot reach $v$ ? This is the case where $j=i+1$. 
Thus, in this case, the boundary interval $[c_{2i},c_{2i+2}]$ is inward, hence 
Lemma \ref{lem:boundaryinterval} tells us that there must exist an inward corner 
incident to this interval, which has to be $c_{2i+1}$. \emph{Two situations then
arise}: either this inward corner can reach $v$, or it cannot, as illustrated in Figure \ref{fig:coins_adjacents}. If it cannot, a
bridge blocks its path, and this creates a $1$-alternating hypermap that can
be detached from our overall hypermap, and this is a bijective operation as explained in the 
last paragraph. But if $c_{2i+1}$ can already reach $v$, then there is no detached hypermap and this will be associated with a weight $1$.
Overall, the weight of the detachable hypermap is either $1$, if $c_{2i+1}$ can already reach $v$, or that of a
$1$-alternating hypermap if it cannot. The total associated weight is therefore the sum
of the weights of these two disjoint cases, namely $H_1(i+1\to i+1)=H_1(x_{i+1},y_{i+1})$.
This is precisely why we added a $\delta_{k,1}$ in the definition of $H_k$,
see Equation \eqref{eq:Hkdef}.
\\

Our algorithm thus allows us to dissect bijectively a pointed $k$-alternating hypermap
into several sub-hypermaps. The first is the component containing the marked
vertex; it now only contains marked corners (inward and outward) that can reach $v$,
according to Lemma \ref{lem:accessible}. It is therefore an
accessibly pointed alternating hypermap.
\\
At each of the inward corners of this accessibly pointed hypermap, we can attach a
bridge oriented outward, together with a general alternating hypermap
counted by an $H_i$. The weight of the bridge is thus accounted for by the
fact that we added an extra term $\prod_{i=1}^k x_i^{-1}y_i^{-1}$ in the
definition of our hypermap generating series. These terms thus allow us to
glue and unglue various hypermaps while simplifying our formulas. Figure \ref{fig:preuve_totale} 
gives a sketch of the global picture of this decomposition.
\\

It only remains to sum over $p$ and to complete the expression with
the weights of the boundary edges, read counterclockwise, to recover
the formula of Theorem \ref{th:dHk/dt}, which we restate here:
$$
\boxed{\dpart{H(1\to k)}{t}=\sum_{p=1}^k\sum_{1\leqslant i_1<\dots<i_p\leqslant k}A_p(x_{i_1},\dots,x_{i_p};y_{i_1-1},\dots,y_{i_p-1})\prod_{j=1}^{p}H(i_j\to i_{j+1}-1),}
$$
where we took $i_{p+1}:=i_1$ for the sake of keeping the invariance by rotation.

\begin{rem}\label{rem:monochromatic}
Dealing with the case of a monochromatic boundary is not hard at all. Indeed, if the boundary 
is monochromatic, then there cannot exist any bridge, hence the pointed hypermap must be accessibly pointed. 
This is a special case of our theorem and proof as it can be seen as a particular case of the subcase $p=k$ by merging all the corners to $c_0$ and
 erasing all the satellite sub-hypermaps that are attached to 
the accessible component, which corresponds to taking a global weight $1$ for the product $\prod_{j=1}^{p}H(i_j\to i_{j+1}-1)=\prod_{j=1}^k H_1(x_j;y_j)$.
\end{rem}

\subsection{Further notions to integrate Equation \eqref{eq:dHk/dt}}\label{sec:hk}

The goal of this subsection is now, starting from Theorem
\ref{th:dHk/dt}, to prove Corollary \ref{cor:hk} by finding a nice 
way to integrate with respect to the vertex weight $t$. To do this, we shall
begin by proving a lemma allowing us to simplify the writing of the
function $A_p$ by expressing it in terms of the $A_1$'s. We will then use
strong induction to prove the result, taking the opportunity to give a way
to recursively compute the coefficients $C_{\sigma}$.

\subsubsection{Linearization of $A_p$}\label{sec:ap_linearise}

We recalled in Section \ref{sec:resultBEL} Relation \eqref{eq:ak}:
$$
\forall p\in\N^*, A_p(x_1,\ldots,x_p;y_1,\ldots,y_p) =[z^0] \frac{1}{\prod_{i=1}^p (x_i-x(z))(y_i-y(z))},
$$
and thus in particular the relation for $p=1$:
$$
A_1(x_1;y_1) =[z^0] \frac{1}{(x_1-x(z))(y_1-y(z))}.
$$
We shall now show that these two relations, the one for $A_p, p\in\N^*,$ and
the one for $A_1$, are intrinsically related combinatorially. I take this
opportunity to thank Gilles Schaeffer, who opened my eyes to this property,
originally appearing in \cite{BEL26} but only via an algebraic partial
fraction decomposition computation. Note then that this section is optional for the understanding of the whole proof of Corollary \ref{cor:hk} from Theorem \ref{th:dHk/dt}.
The only thing that will really matter is the Equation \eqref{eq:ap_linearise}
\\

\begin{figure}
    \centering
    \includegraphics[width=\textwidth]{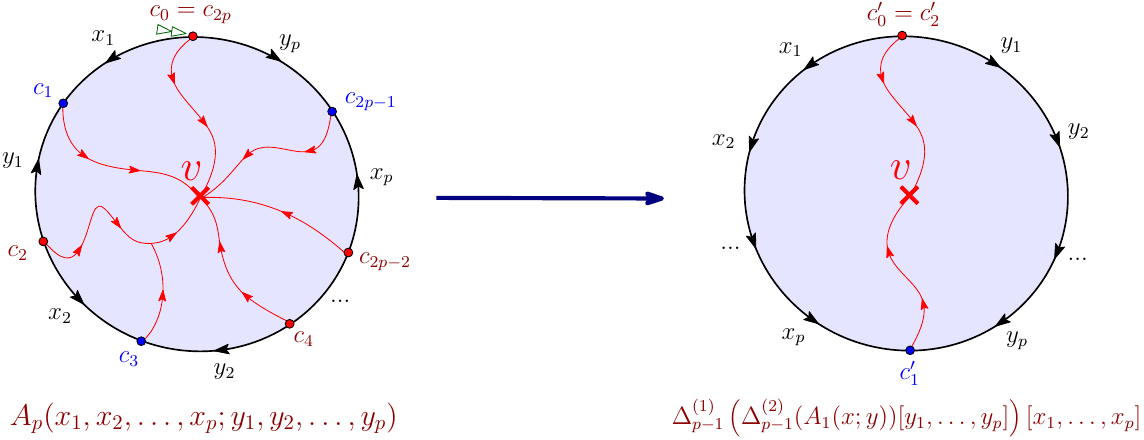}
    \caption{Sketch of the decomposition of a $p$-alternating accessibly pointed hypermap to a $1$-alternating one with multiple weights.}
    \label{fig:ap_to_a1}
\end{figure}

Let us return to the definition of $A_p$: it is the generating series of
accessibly pointed $p$-alternating hypermaps. We recall to the interested reader that this type of hypermap can be
decomposed according to sequences of slices of various weights, namely the
weights of the edges appearing on the boundary of the hypermap, which leads to Formula \eqref{eq:ak} as proved in \cite{BEL26}. This
decomposition is then symmetric in the $x_i$'s and the $y_i$'s. We can thus write it in the form:
\begin{equation}\label{eq:ak_prodxi*prodyi}
\forall p\in\N^*, A_p(x_1,\ldots,x_p;y_1,\ldots,y_p) =[z^0] \frac{1}{\prod_{i=1}^p (x_i-x(z))}\frac{1}{\prod_{i=1}^p(y_i-y(z))}.
\end{equation}
We have thus only grouped the $x_i$'s together and the $y_i$'s together.
This formula thus translates into the presence of an accessible
$1$-alternating hypermap, but whose black boundary interval weights are not
identical and are counted by the $x_i$'s, while those of the white boundary
interval would be counted by the $y_i$'s, see Figure \ref{fig:ap_to_a1} for
more clarity. We also refer the reader to \cite{BOUTTIER_2014,AB2025,BEL26}, which allows for
a better understanding of the cutting and gluing of slices, providing the
explicit bijective argument that makes our approach rigorous.

In order to translate this remark into an equation, we must go further. We
now introduce several operators $\Delta_n, n\in\N^*$ acting linearly on the
ring of formal series in one variable with coefficients in $\mathcal{R}$,
such that $\Delta_n$ takes as input a series in
$\frac{1}{a}\mathcal{R}[\![a^{-1}]\!]$ and takes values in
$\frac{1}{a_1\cdots a_{n+1}}\mathcal{R}[\![a_1^{-1},\dots,a_{n+1}^{-1}]\!]$, where $a$ and the $a_i$s are some formal variables.
We define them recursively as follows, for
$f\in\frac{1}{a}\mathcal{R}[\![a^{-1}]\!]$:
\begin{equation}\label{eq:def_deltan}\left\{\begin{array}{lll}
&\Delta_1(f)[a_1,a_2]:=\displaystyle\frac{1}{a_1a_2}\frac{f(a_1)-f(a_2)}{a_1^{-1}-a_2^{-1}}\,\,\,\text{and, for all $n\geqslant1$,}\\
&\\
&\Delta_n(f)[a_1,\dots,a_{n+1}]:=\Delta_1\big(\Delta_{n-1}(f)[a_1,\dots,a_{n-1},\bullet]\big)[a_n,a_{n+1}]\\
&\\
&\hspace{3.485cm}=\displaystyle\frac{1}{a_na_{n+1}}\frac{\Delta_{n-1}(f)[a_1,\dots,a_{n-1},a_{n}]-\Delta_{n-1}(f)[a_1,\dots,a_{n-1},a_{n+1}]}{a_n^{-1}-a_{n+1}^{-1}}.
\end{array}\right.
\end{equation}
The usefulness of this generalized divided difference operator is to
describe a single combinatorial object possessing multiple weights to be
counted. Indeed, we have the equality:
$$
\Delta_1\left(f(a):=\sum_{n=0}^{+\infty}\frac{f_n}{a^{n+1}}\right)[a_1,a_2]=\frac{1}{a_1a_2}\sum_{n=0}^{+\infty}\left(\frac{1}{a_1^{n}}+\frac{1}{a_1^{n-1}a_2}+\frac{1}{a_1^{n-2}a_2^2}+\dots+\frac{1}{a_2^{n}}\right)\cdot f_n.
$$
We recover the same series as $f(a)$, but where the weight $a^{-n}$ in front
of the term $f_n$ is split in two: one part will have a weight
$a_1^{-1}$ and the other a weight $a_2^{-1}$, such that there are $n$
weights in total (the extra $\frac{1}{a_1a_2}$ is forgotten here). Recall that we must add an additional weight
$a_1^{-1}a_2^{-1}$ in order to follow the conventions we have taken from
the beginning, see Equation \eqref{eq:Hkdef}. But this is exactly the
operator that corresponds to our description of Formula
\eqref{eq:ak_prodxi*prodyi}! This idea generalizes to the interpretation of
$\Delta_n(f(a))[a_1,\dots,a_n]$ from Equation \eqref{eq:def_deltan} as counting the objects described by $f$ no
longer by a single weight $a^{-1}$ but by several weights
$a_1^{-1},\dots,a_{n+1}^{-1}$. This is thus exactly the object we need in
order to understand $A_p$, in view of what we have stated previously.
\\

Thus, in what follows, our series $f(a)$ is replaced by $A_1(x;y)$. We will
moreover define the induced operators $\Delta_n^{(1)}$ and $\Delta_n^{(2)}$,
which will act respectively on the variable $x$ and on the variable $y$ of the
input series. We thus have the relation:
$$
A_2(x_1,x_2;y_1,y_2)=\Delta_1^{(1)}\left(\Delta_1^{(2)}(A_1(x;y))[y_1,y_2]\right)[x_1,x_2].
$$
This gives, decomposing step by step:
$$
\Delta_1^{(2)}(A_1(x;y))[y_1,y_2]=\frac{1}{y_1y_2}\frac{A_1(x;y_1)-A_1(x;y_2)}{y_1^{-1}-y_2^{-1}},
$$
then
\begin{align*}
\Delta_1^{(1)}\left(\Delta_1^{(2)}(A_1(x;y))[y_1,y_2]\right)[x_1,x_2]&=\frac{1}{y_1y_2}\frac{\Delta_1^{(1)}(A_1(x;y_1))[x_1,x_2]-\Delta_1^{(1)}(A_1(x;y_2))[x_1,x_2]}{y_1^{-1}-y_2^{-1}}\\
&=\frac{1}{y_1y_2}\frac{1}{y_1^{-1}-y_2^{-1}}\Bigg(\frac{1}{x_1x_2}\frac{A_1(x_1;y_1)-A_1(x_2;y_1)}{x_1^{-1}-x_2^{-1}}\\
&\hspace{4cm}-\frac{1}{x_1x_2}\frac{A_1(x_1;y_2)-A_1(x_2;y_2)}{x_1^{-1}-x_2^{-1}}\Bigg)\\
&=\frac{A_1(x_1;y_1)+A_1(x_2;y_2)-A_1(x_1;y_2)-A_1(x_2;y_1)}{(x_1-x_2)(y_1-y_2)}.
\end{align*}
And it is exactly this formula which, in \cite{BEL26}, allowed us to
recover the expression of $H_2$ from Equation \eqref{eq:h2} starting from
Equation \eqref{eq:dH2/dt} by integrating!
\\

Now that the notations are understood, we can write the generalization of
our formula for an arbitrary $p\in\N^*$:
$$
A_p(x_1,\dots,x_p;y_1,\dots,y_p)=\Delta^{(1)}_{p-1}\left(\Delta^{(2)}_{p-1}(A_1(x;y))[y_1,\dots,y_p]\right)[x_1,\dots,x_p].
$$

More explicitly, we can express each term of the final linearization in the
form
\begin{equation}\label{eq:ap_linearise}
\boxed{A_p(x_1,\dots,x_p;y_1,\dots,y_p)=\sum_{1\leqslant i,j\leqslant p}\frac{A_1(x_i;y_j)}{\prod_{a\neq i}(x_a-x_i)\prod_{b\neq j}(y_b-y_j)},}
\end{equation}
which is the formula that will be used in what follows. This formula
follows from a quick induction using the recursive definition of the
$\Delta_n$'s that is not worth developing here.

\subsubsection{A word on plane permutations}\label{sec:planaire}

As stated in Corollary \ref{cor:hk}, we will need the notion of plane
permutation before beginning the proof of this corollary.
\begin{Def}\label{def:planaire}
Let $\sigma\in\mathfrak{S}_k$ be a permutation of $[\![1,k]\!]$. Let us also
consider the graph $G_{\sigma}$ whose vertices are those of a regular
$2k$-gon, labeled alternately $x_1, y_1, x_2, y_2, \dots, x_k, y_k$
in a counterclockwise direction, and whose edges are the segments $[x_i, y_{\sigma(i)}]$
for every $i \in [\![1, k]\!]$, and the $2k$ sides of the $2k$-gon.
\\

We then say that $\sigma$ is plane if, and only if, the graph $G_{\sigma}$
is plane. In other words, if none of its edges cross, which then forms a
plane map. We will denote by $\mathfrak{P}_k\subseteq\mathfrak{S}_k$ their
set.

\begin{figure}[h!]
    \centering
    \begin{subfigure}[b]{0.4\textwidth}
        \includegraphics[width=\textwidth]{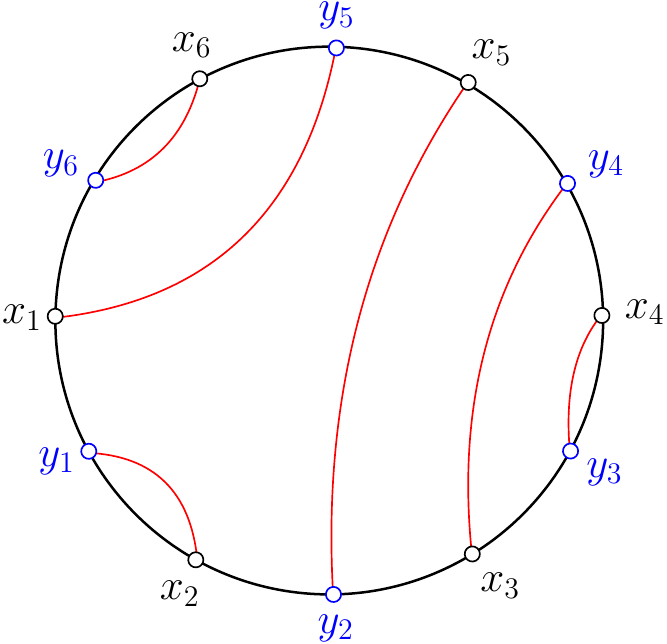}
        \caption{The permutation $(514326)$ is plane,\\ associated with the word (()((())))().}
    \end{subfigure}
    \hspace{2cm}
    \begin{subfigure}[b]{0.4\textwidth}
        \includegraphics[width=\textwidth]{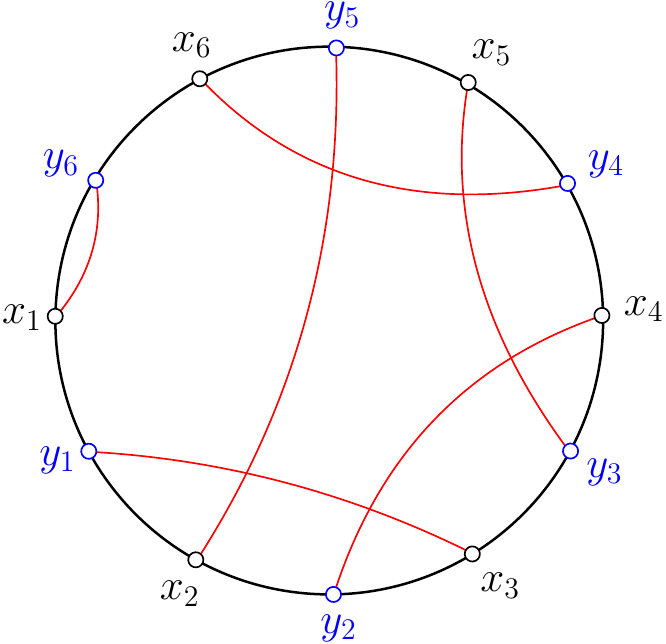}
        \caption{The permutation $(651234)$ is not plane.}
    \end{subfigure}
    \caption{Graphs associated with two permutations of size $k=6$. Here the edges are curved for an aesthetic purpose.}
\end{figure}
\end{Def}

\begin{rem}
In the literature, one often translates such a structure in terms of \emph{non-crossing chord diagrams}. The reader familiar with this notion can keep in mind that a plane 
permutation of order $k$ is nothing but a function that can be represented as a non-crossing chord diagram with $k$ chords.
\end{rem}

\begin{rem}\label{rem:planaire_generale}
In what follows, we will generalize this definition to any finite subset
$E\neq\emptyset$ of $\N^*$, rather than just $[\![1,k]\!]$. To do this, one can send $E$
bijectively, respecting the order of its elements, onto $[\![1,k]\!]$, where
$k:=\vert E\vert$, and this naturally generalizes our definition of a
plane permutation. We will thus denote by $\mathfrak{P}_E$ the set of
plane permutations on $E$.
\end{rem}
Let us now see a quick characterization of plane permutations in order to
get a feel for this notion. Recall that a \emph{balanced string of brackets} is a word
made up of the two symbols '(' and ')' such that no bracket can be
closed before being opened. In other words, when reading the word, one must
always have counted more '(' than ')' and the total number of each of these symbols should be equal.
\begin{prop}\label{prop:plan_parenthese}
The set of plane permutations of $[\![1,k]\!]$ is in bijection with the
set of balanced string of brackets of size $2k$. In particular,
$\vert\mathfrak{P}_k\vert=\displaystyle Cat(k):=\frac{\binom{2k}{k}}{k+1}$.
\end{prop}

\begin{proof}
    Start from a plane permutation $\sigma$ and its associated
    plane graph. Also consider the list of characters
    $[x_1,y_1,\dots,x_k,y_k]$, whose elements will be modified afterward as we
    read the graph. \\

    Let us begin our induction, with $x_1$, which is connected to
    $y_{\sigma(1)}$. We then change $x_1$ into '(' in the list and
    $y_{\sigma(1)}$ into ')'. We will then say that $x_1$ and $y_{\sigma(1)}$
    have already been visited.
    \\

    We then continue our circular reading, counterclockwise, of the graph
    at the next unvisited site. Say it is $y_{1}$. We then replace
    $y_1$ with '(' in our list and $x_{\sigma^{-1}(1)}$ with ')'. We
    continue going around our graph in this way until we have replaced all
    the elements of the list with brackets. Moreover, the final list
    forms a balanced string of brackets by construction, since we have always
    started with an opening bracket before closing it as we traversed
    the list.\\

    The construction of the inverse of our function is done by exactly
    retracing the reasoning from the bottom up, looking for when each
    bracket opens and closes in order to construct $\sigma$, which
    concludes the proof.

\end{proof}

Of course, this result generalizes to any subset $E$ of $\N^*$ of size $k$
instead of $[\![1,k]\!]$, by Remark \ref{rem:planaire_generale}. It only
remains for us to prove the key result of this section, which will be used
in the proof of Corollary \ref{cor:hk}.

\begin{lem}\label{lem:permutation}
Let $I$ be a finite subset of $[\![1,k]\!]$ and $p\geqslant2$ its
cardinality; we will denote by $I[j]$ the $j$-th element of $I$ considered
in strictly increasing order. Let now $P_1,\dots, P_p$ be the subsets of
$[\![1,k]\!]$ such that $P_j:=[\![I[j],I[j+1]-1]\!], j\in[\![1,p-1]\!]$ and
$P_p:=[\![I[p],k]\!]\sqcup[\![1,I[1]-1]\!]$.\\

Let us finally consider $\sigma_1,\dots,\sigma_p$ plane permutations such
that $\sigma_j\in\mathfrak{P}_{P_j}$ for every $1\leqslant j\leqslant p$.
\\
Then the function $\sigma_I\colon [\![1,k]\!]\to [\![1,k]\!]$ defined by
$\sigma_I(i):=\sigma_j(i)$ for every $i\in P_j, j\in[\![1,p]\!]$, is a
plane permutation.
\end{lem}

\begin{proof}
    It is clear, by construction, that $\sigma_I$ is injective and
    surjective, since the $\sigma_i$'s are and $I=\bigsqcup_{i=1}^pP_i$. It is thus indeed a permutation.
    \\

    By rotational invariance of the notion of planarity (considering the
    graph associated with a plane permutation), we may assume that
    $I[1]=1$ without loss of generality. Thus $P_1$ starts at $1$ and $P_p$ ends at $k$.
    \\

    By Proposition \ref{prop:plan_parenthese}, we can associate to each
    $\sigma_i$ a balanced string of brackets. By construction of $\sigma_I$, it
    suffices to concatenate all these words to form the string of brackets
    associated with $\sigma_I$. This word is then balanced, since
    all its sub-words are, and thus $\sigma_I$ is plane, which concludes
    the proof.
\end{proof}

\begin{figure}[h!]
    \centering
    \includegraphics[width=0.8\textwidth]{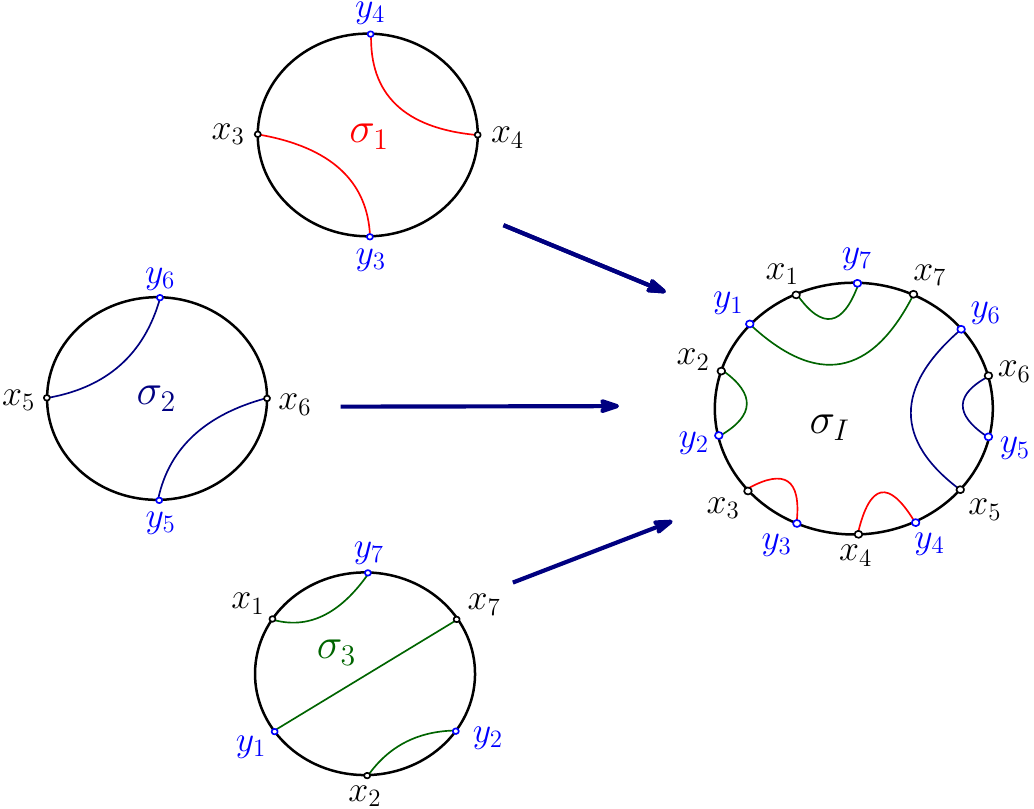}
    \caption{Construction of the graph of $\sigma_I, I=\{3,5,7\}\subseteq[\![1,7]\!]$ from the three permutations $\sigma_1:=(3,4)\mapsto (3,4),\sigma_2:=(5,6)\mapsto (6,5)$ and $\sigma_3=(7,1,2)\mapsto (1,7,2)$.}
    \label{fig:concatenation_permut}
\end{figure}

Informally speaking, this construction of $\sigma_I$ simply consists of a
"union" of the plane graphs associated with a tuple of plane
permutations. We have thus constructed a large permutation $\sigma_I$ from
plane permutations on $P_1,\dots, P_p$, and this has formed for us a
plane permutation on $[\![1,k]\!]$. Figure \ref{fig:concatenation_permut} gives a sketch of this idea.

In what follows, we will denote by $\mathfrak{P}_{I,k}$ the set of plane
permutations that can be decomposed as such $\sigma_I$, considering $I$ as a subset of $[\![1,k]\!]$. \\
Now we need to define 
a specific concept about plane permutations that will be very useful in the next section.
\begin{Def}
    Let $\emptyset \neq E\subseteq \N^*$ be a finite set, $\emptyset\neq A\subseteq E$ a subset of $E$ and $\sigma\in\mathfrak{P}_E$ be a plane permutation on $E$. 
    We say that $A$ is \emph{stable} under $\sigma$ if $\sigma(A)=A$.
    \\
    In this case, Definition \ref{def:planaire} implies that the restriction $\sigma\colon A\to A$ is a plane permutation on $A$. 
\end{Def}

One might wonder from Lemma \ref{lem:permutation} whether every plane permutation can be associated with
some subset $I$ of $[\![1,k]\!]$, in other words whether the equality
$\displaystyle\bigcup_{I\subseteq [\![1,k]\!]}\mathfrak{P}_{I,k}\overset{?}{=}\mathfrak{P}_k$ holds. It is exactly what the following proposition is about. 
\begin{prop}\label{prop:permutation}
Let $k\geqslant2$ be an integer and $\gamma_k\in\mathfrak{S}_k$ be the cycle $(k,\,k-1,\dots,1)$.
\\
Then we have the set equality:
$$
\displaystyle\bigcup_{I\subseteq [\![1,k]\!]}\mathfrak{P}_{I,k}=\mathfrak{P}_k\backslash\{\gamma_k\},
$$
where the union on the left-hand side is \underline{not} disjoint.
\end{prop}

\begin{proof}
Let us consider $\sigma\in\mathfrak{P}_k$ such that $\sigma\neq\gamma_k$. Hence there exists $i$ such that $\sigma(i)\neq i-1$ (where $0\equiv k$ by convention).
Now let us define $I:=\{i,\sigma(i)+1\}$ (which hence contains at least two elements) where we assumed that $I[1]=i$ and $I[2]=\sigma(i)+1\leqslant k$ without loss of generalities.
\\

We can then define $P_1:=[\![I[1],I[2]-1]\!]=[\![i,\sigma(i)]\!]$ and $P_2:=[\![I[2],k]\!]\cup[\![1,I[1]-1]\!]$ which is $[\![\sigma(i)+1,k]\!]\cup[\![1,i-1]\!]$. Then, since $\sigma$ is plane, there cannot exist any $j$ in $P_1$ that can be sent to $P_2$ because the 
edge $x_j\sim y_{\sigma(j)}$ would cross the edge $x_i\sim y_{\sigma(i)}$ in $G_{\sigma}$. Hence $P_1$ is stable under $\sigma$ and so is $P_2$ since $\sigma$ is bijective. Then $\sigma=\sigma_I$, where $\sigma_1:=\sigma_{|P_1}$ and $\sigma_2:=\sigma_{|P_2}$, so $\sigma\in\mathfrak{P}_{I,k}$, which gives the converse inclusion.
\\

Finally, $\gamma_k$ is a plane permutation that makes no strict subset of $[\![1,k]\!]$ stable under itself, hence it cannot belong to any $\mathfrak{P}_{I,k}$. This concludes the proof as the inclusion of all $\mathfrak{P}_{I,k}$ into $\mathfrak{P}_k$ is given by Lemma \ref{lem:permutation}.
\end{proof}

This is why the cycle $\gamma_k$ will have a very important role in the next section.

\subsection{An algebraic formula for $H_k$ as a function of $H_1$}\label{sec:2.4}

In this section we begin the proof by strong induction of Corollary
\ref{cor:hk}. We will start by stating the strong induction property to be
proven, arising from Proposition \ref{prop:recu}, which directly implies
Corollary \ref{cor:hk}. We will then support Theorem \ref{th:dHk/dt} with
some computations in order to prove Lemma \ref{lem:introcalphabeta}.
Finally, the third subsection will aim to finalize the inductive step,
which will conclude our proof.

\subsubsection{Introducing the induction property}\label{sec:intro_csigma}

Let us begin with a few definitions. For each integer $n\in\N^*$, there
exists a function $C$ taking as input an element of $\mathfrak{P}_n$ and
with values in the ring of rational fractions in the variables
$x_1,y_1,\dots,x_n,y_n$ with integer coefficients, defined by the formulas
\begin{equation}\label{eq:recu_csigma}
\forall \sigma\in\mathfrak{P}_n, \left\{\begin{array}{lll}
    &C_{\sigma}(1\to1)=1\,\,\,\text{if $n=1$},\\
    &\\
& \hspace{-0.75cm}\left.\begin{array}{lll}
&C_{\sigma}(1\to n)=\displaystyle\sum_{j=1}^{n-1}\frac{C_{\sigma}(1\to j)\,C_{\sigma}(j+1\to n)}{(x_{j+1}-x_{1})(y_n-y_{\sigma(1)})}\,\,\,\text{if $\sigma(1)\neq n$},\\
&C_{\sigma}(x_1,\dots,x_n;y_1,\dots,y_n)=C_{\gamma_n^{-1}\circ\sigma\circ\gamma_n}(x_2,\dots,x_n,x_1;y_2,\dots,y_n,y_1),\\
&\\
&\displaystyle\sum_{\sigma\in\mathfrak{P}_n}C_{\sigma}(1\to n)=0.
\end{array}
\right\}\,\,\,\text{if $n>1$,}
\end{array}
\right.
\end{equation}
where we have defined, for $1\leqslant a\leqslant b\leqslant n$ and
$\sigma\in\mathfrak{P}_n$, $C_{\sigma}(a\to b):=C_{\sigma_{\mid[\![a,b]\!]}}(a\to b)$
where $\sigma_{\mid[\![a,b]\!]}$ is the restriction of $\sigma$ from
$[\![a,b]\!]$ to $[\![a,b]\!]$ if $[\![a,b]\!]$ is stable under $\sigma$, otherwise we set $C_{\sigma}(a\to b):=0$. Also, we recall the convention 
of Theorem \ref{th:dHk/dt} which also holds for $C$ and leads to the definition 
$C_{\sigma}(1\to n):=C_{\sigma}(x_1,\dots,x_n;y_1,\dots,y_n)$.
\\

The first two equations allow us to define $C$ by induction on the
size of the permutations, if $1$ is not sent to $n$. The third one allows
us to circumvent the condition $\sigma(1)\neq n$ by rotational invariance,
letting us reduce to any $i$ such that $\sigma(i)\neq i-1$, a condition
which then guarantees the existence of intervals of integers stable under
$\sigma$ strictly contained in $[\![1,n]\!]$ as explained in Proposition \ref{prop:permutation}.
 Thus it only remains to define $C_{\sigma}$ for the cycle $\sigma=\gamma_n=(n,\,n-1,\dots,1)$ in order to
fully determine $C$, which is exactly the role of the last equation. Note that we only used the third equation for $\sigma(1)=n$ 
so one might ask if $C$ is well-defined, and indeed it is no more than a redundancy but the proof is none of our interest here so we refer the curious reader to \cite{EO2005} for further details.

\begin{rem}
$\bullet$ All these relations defining $C$ are not exactly the ones that already
appeared and were shown analytically in \cite{EO2005}, some of which
were taken up again in \cite{Eynard2016}. The difference is simply a
convention: their white faces are our black faces and vice versa, and
moreover they turn clockwise while we turn counterclockwise. To go from
one to the other, one must swap $\sigma$ with $\sigma^{-1}$, $x_i$ and
$y_i$, and finally flip our diagrams.

$\bullet$ Another difference here is that our method, originally based on a
combinatorial reasoning, is purely formal. It in particular allows us to
show the strength of the notion of accessibility for proving powerful
theorems.

$\bullet$ Let us take this opportunity to point out a small typo appearing on page
401 of the book \cite{Eynard2016}. Indeed, the relation given for the $C_{\sigma}$'s cannot hold if
$\sigma^{-1}(n)=1$, which in particular gives the problem associated with
the cycle $\gamma_n$ described previously. However the first article stating this property, which is 
\cite{EO2005}, states it faithfully.
\end{rem}

In fact, these relations defining $C$ are a bit tricky to understand. Let us show a useful lemma 
that will be useful later on in our future computation, and which make easier the computation of $C_{\sigma}$ for a fixed planar permutation $\sigma$, just by taking a quick look on its associated graph $G_{\sigma}$.

\begin{lem}\label{lem:csigma_face12}
Let $\sigma\in\mathfrak{P}_k$ be a plane permutation such that $\sigma(1)\neq1,k$. Then
$$
C_{\sigma}(1\to k)=C_{\sigma}(1\to \sigma(1))\cdot C_{\sigma}(x_1,y_{\sigma(1)},x_{\sigma(1)+1},y_{\sigma(1)+1},\dots,x_k,y_k).
$$
Using our earlier notation, this formula generalizes to
\begin{equation}\label{eq:csigma_face12}
C_{\sigma}(F)=C_{\sigma}(F_1)\,C_{\sigma}(F_2),
\end{equation}
where $F_1$ and $F_2$ are, respectively, the two faces obtained by cutting the face $F$ along an edge $e:=(i,\sigma(i))\notin\{(i,i),(i,i-1)\}$, when such an edge exists.
\end{lem}

\begin{proof}
We prove the first identity using the recursive Relation \eqref{eq:recu_csigma} defining the $C_{\sigma}$; the second follows easily from the first by labeling the $x_i$ and $y_i$ found around our two faces and using rotational invariance of the $C_{\sigma}$ to start wherever we please.
\\

We write, isolating the $j=1$ term of the first sum and noting that the $y_{\sigma(1)}$ in Equation \eqref{eq:recu_csigma} remains fixed — since it corresponds to $y_i$ where $i$ is the image under $\sigma$ of the index of the first $x$ appearing in $C_{\sigma}$, namely $1$ here:
\begin{align*}
C_{\sigma}(1\to \sigma(1))\cdot C_{\sigma}(x_1;y_{\sigma(1)},x_{\sigma(1)+1},\dots,y_k)&= C_{\sigma}(1\to \sigma(1))\cdot \Bigg(\frac{C_{\sigma}(x_1,y_{\sigma(1)})\,C_{\sigma}(\sigma(1)+1\to k)}{(x_{\sigma(1)+1}-x_1)(y_k-y_{\sigma(1)})}\\
&+\sum_{j=\sigma(1)+1}^{k-1}\frac{C_{\sigma}(x_1,y_{\sigma(1)},\dots,y_{j})\,C_{\sigma}(j+1\to k)}{(x_{j+1}-x_1)(y_k-y_{\sigma(1)})}\Bigg).
\end{align*}
By strong induction (the base case being immediate), we may form $C_{\sigma}(1\to j)$ by merging $C_{\sigma}(1\to\sigma(1))$ and $C_{\sigma}(x_1,y_{\sigma(1)},\dots,y_{j})$. It remains only to note that $C_{\sigma}(x_1;y_{\sigma(1)})=1$, and to gather everything under the same sum (the first term corresponding to $j=\sigma(1)$), which gives:
$$
C_{\sigma}(1\to \sigma(1))\cdot C_{\sigma}(x_1,y_{\sigma(1)},x_{\sigma(1)+1},\dots,y_k)=\sum_{j=\sigma(1)}^{k-1}\frac{C_{\sigma}(1\to j)\,C_{\sigma}(j+1\to k)}{(x_{j+1}-x_1)(y_k-y_{\sigma(1)})}=C_{\sigma}(1\to k),
$$
the last equality following from the fact that the intervals $[\![1,j]\!], j\in[\![1,\sigma(1)-1]\!]$, are not stable under $\sigma$ (since the image of $1$ under $\sigma$ lies outside them), so that the corresponding $C_{\sigma}(1\to j)$ all vanish as stated after the Equations \eqref{eq:recu_csigma} defining $C$. This concludes the proof.
\end{proof}
\begin{rem}
This lemma simplifies the computation of the coefficients $C_{\sigma}$. Indeed, consider the plane graph associated with $\sigma$: the lemma tells us that $C_{\sigma}$ can be computed by reducing to the restrictions of $\sigma$ on either side of the edge $x_1\sim y_{\sigma(1)}$, whenever $\sigma(1)\neq 1,k$. Repeating this process for every $i$ such that $\sigma(i)\neq i,i-1$, one reduces the computation of $C_{\sigma}$ to its values on the faces of the graph $G_{\sigma}$. Now, the restriction of $\sigma$ to each such face necessarily sends an $x_i$ to one of its neighbors. Since this restriction is itself a permutation, it must therefore send every $x_i$ either to the $y_j$ immediately to its right, or to its left neighbor. In other words, the graph of the restriction coincides either with that of the identity or with that of some $\gamma_p$.
\\
Since the graph of $\gamma_p$ is the same as that of the identity, only traversed in the opposite direction, the symmetry of $H_k$ under this flip allows us to reduce the knowledge of $C_{\sigma}$ entirely to that of $C_{\mathrm{Id}}$, evaluated on each face of its associated plane graph and read in one direction or the other depending on whether we are dealing with $\gamma_p$ or $\mathrm{Id}$.
\\
This decomposition is in fact the one obtained by Eynard and Orantin in \cite{EO2005}. They already expressed $C_{\sigma}$ in terms of what they call "face functions $F_p$", and this remains to this day the simplest known way of determining the $C_{\sigma}$. However, it will not be needed here, so we do not elaborate further; we refer the interested reader to that article for more details.
\end{rem}

Let us now observe what role these coefficients play in our problem. We will
define, for every permutation $\sigma\in\mathfrak{P}_n, n\in\N^*$ and
$\alpha,\beta\in[\![1,n]\!]$:
\begin{equation}\label{eq:calphabeta}
C_{\sigma}^{\alpha,\beta}(1\to n):=\sum_{\underset{\vert I\vert\geqslant2;\, \alpha,\beta+1\in I}{I\subseteq[\![1,n]\!]/}}\frac{\prod_{j=1}^{\vert I\vert}C_{\sigma}(I[j]\to I[j+1]-1)}{\displaystyle\prod_{\underset{i\neq\alpha}{i\in I}}(x_{i}-x_{\alpha})\prod_{\underset{j\neq\beta+1}{j\in I}}(y_{j-1}-y_{\beta})}
\end{equation}
where the elements of $I$ are considered modulo $n$. In particular, one notices that the coefficient
$C_{\sigma}^{1,n}(1\to n)$ takes a particular form that will be useful to
us later:
\begin{equation}\label{eq:c1ndef}
C_{\sigma}^{1,n}(1\to n):=\sum_{\underset{\vert I\vert\geqslant2;\, 1\in I}{I\subseteq[\![1,n]\!]/}}\frac{C_{\sigma}(1\to I[2]-1)\cdots C_{\sigma}(I[\,\vert I\vert\,]\to n)}{\displaystyle\prod_{\underset{i\neq1}{i\in I}}(x_{i}-x_{1})\prod_{\underset{j\neq1}{j\in I}}(y_{j-1}-y_{n})},
\end{equation}
using $C_{\sigma}(n+1\to1)=C_{\sigma}(1\to1)=1$. Since the variables $x_i$
and $y_i$ are formal, we deduce the following definition for
$1\leqslant \alpha\leqslant \beta \leqslant n$ and
$\sigma\in\mathfrak{P}_{[\![\alpha,\beta]\!]}$:
\begin{equation}\label{eq:c1k}
C_{\sigma}^{\alpha,\beta}(\alpha\to \beta):=\sum_{\underset{\vert I\vert\geqslant2;\, \alpha\in I}{I\subseteq[\![\alpha,\beta]\!]/}}\frac{C_{\sigma}(\alpha\to I[2]-1)\cdots C_{\sigma}(I[\,\vert I\vert\,]\to \beta)}{\displaystyle\prod_{\underset{i\neq\alpha}{i\in I}}(x_{i}-x_{\alpha})\prod_{\underset{j\neq\beta+1}{j\in I}}(y_{j-1}-y_{\beta})}.
\end{equation}

In fact the reason we introduced the
$C_{\sigma}^{\alpha,\beta}$ is because they are strongly related to the $C_{\sigma}$, which themselves are related to the generating function $H_k$. And 
they will have a particular role when we will aim to make calculations on the right-hand side of Theorem \ref{th:dHk/dt}, as Lemma 
\ref{lem:introcalphabeta} illustrates. 
\\

Let us now state the fundamental proposition of this section which explains exactly the relation between the $C_{\sigma}^{\alpha,\beta}$'s and the $C_{\sigma}$'s. 
\begin{prop}\label{prop:recu}
For every $k\in\N^*$, $\sigma\in\mathfrak{P}_k$ and $\alpha,\beta\in[\![1,k]\!]$,
we can compute the coefficients $C_{\sigma}^{\alpha,\beta}$:
\begin{equation}
C_{\sigma}^{\alpha,\beta}(1\to k)=\left\{\begin{array}{lll}
    &C_{\sigma}(1\to k)\,\,\,\text{if $\beta=\sigma(\alpha)$},\\
    &-\,C_{\sigma}(1\to k)\,\,\,\text{if $\beta+1=\alpha$},\\
    &0\,\,\,\text{otherwise}.
\end{array}
\right.
\end{equation}
Finally, these functions $C_{\sigma}$ are directly related to $H_k$ by the
relation
\begin{equation}\label{eq:hk}
H(1\to k)=\sum_{\sigma\in\mathfrak{P}_k}C_{\sigma}(1\to k) \prod_{i=1}^kH_1(x_i;y_{\sigma(i)}).
\end{equation}
\end{prop}

It is clear that this Proposition directly implies Corollary \ref{cor:hk}. We will prove it by \textbf{induction} on $k$.

\paragraph{Base case} We can begin by quickly warming up by proving the
base case, namely $k=1$.\\

In particular, for $k=1$, we have $\mathfrak{P}_1=\{\mathrm{Id}\}$ and we
can only take $\alpha=\beta=1$. We thus have, by Equation \eqref{eq:c1k}
$$
C_{\mathrm{Id}}^{1,1}(1\to1)=\sum_{\underset{\vert I\vert \geqslant 2,1\in I}{I\subseteq[\![1,1]\!]}}\dots=\sum_{I\in \emptyset}\dots=0,
$$
which is exactly what we wanted, since this quantity is equal both to
$C_{\mathrm{Id}}$ and to $-C_{\mathrm{Id}}$, given that $\beta=\mathrm{Id}(\alpha)$
and $\beta+1=\alpha$ with the convention $\beta+1\equiv 1$ if $\beta=n=1$.
\begin{rem}\label{rem:gammak}
Moreover, $C_{\gamma_n}^{\alpha,\beta}=0$ for every $\alpha,\beta\in[\![1,n]\!], n\geqslant1$,
since the only non-empty subset of $[\![1,n]\!]$ stable under $\gamma_n$ is $[\![1,n]\!]$. And this is in fact exactly what we have just
recovered here, since $\gamma_1=\mathrm{Id}$, and is also a special case of Proposition \ref{prop:recu} as the equality $\beta+1=\alpha$ is equivalent to $\beta=\gamma_n(\alpha)$.
\end{rem}

This base case, although necessary, has not allowed us to gain a better
understanding of the definitions of the $C_{\sigma}^{\alpha,\beta}$'s. Let
us therefore move on to the next section, whose goal is to rewrite Theorem
\ref{th:dHk/dt} in the form of Equation \eqref{eq:final_htilde}, thereby
bringing out our $C_{\sigma}^{\alpha,\beta}$'s and making the statement of
Proposition \ref{prop:recu} more explicit. \\

\textbf{We now assume that the induction property defined in Proposition
\ref{prop:recu} holds up to rank} $\mathbf{k-1, k\geqslant2}$. We then seek to prove
it at rank $k$. To do this, we will need to prove a lemma.

\begin{lem}\label{lem:introcalphabeta}
Assuming the property of Proposition \ref{prop:recu} up to rank $k-1$, we
can rewrite Theorem \ref{th:dHk/dt} in the form:
\begin{equation}\label{eq:final_htilde}
\boxed{\displaystyle\dpart{}{t}\left(\frac{H(1\to k)}{\prod_{i=1}^k H_1(x_i;y_{i-1})}\right)=\displaystyle\sum_{\underset{\sigma\neq\gamma_k}{\sigma\in\mathfrak{P}_k}}\displaystyle\frac{\prod_{i=1}^kH_1(x_i;y_{\sigma(i)})}{\prod_{i=1}^k H_1(x_i;y_{i-1})}\displaystyle\sum_{1\leqslant\alpha,\beta\leqslant k}A_1(x_{\alpha};y_{\beta})\cdot C_{\sigma}^{\alpha,\beta}(1\to k).}
\end{equation}
\end{lem}

The impatient reader may refer to Section \ref{sec:conclusion_hk} to see
how this lemma can lead to a proof of our inductive step.

\subsubsection{Proof of Lemma \ref{lem:introcalphabeta}}\label{sec:introcalphabeta}

Let us return to Theorem \ref{th:dHk/dt} and isolate the term
corresponding to $p=1$ in the sum. Setting $y_0:=y_k$, we obtain
$$
\dpart{H(1\to k)}{t}=\sum_{i=1}^kA_1(x_{i};y_{i-1})H(i\to i-1)+\sum_{p=2}^k\dots
$$
and by our convention on the notation $H(i\to i-1)$, this expression
is simply equal to
$$
H(i\to i-1)=H_k(x_i,\dots,x_k,x_1,\dots,x_{i-1};y_i,\dots,y_k,y_1,\dots,y_{i-1})=H(1\to k),
$$
by rotational invariance of the
problem, for all $i\in[\![1,k]\!]$. Using Relation \eqref{eq:dH1/dt}, we get 
$$
\left(\prod_{i=1}^k H_1(x_i;y_{i-1})\right)\,\dpart{}{t}\left(\frac{H(1\to k)}{\prod_{i=1}^k H_1(x_i;y_{i-1})}\right)=\dpart{H(1\to k)}{t}-\left(\sum_{i=1}^n A_1(x_i;y_{i-1})\right)H(1\to k).
$$
Then moving the $p=1$ term to the left-hand side and using the last three equations gives 
$$
\dpart{\tilde{H}(1\to k)}{t}=\frac{1}{\prod_{i=1}^k H_1(x_i;y_{i-1})}\sum_{p=2}^k\sum_{i_1<\dots<i_p}\hspace{-0.2cm}A_p(x_{i_1},\dots,x_{i_p};y_{i_1-1},\dots,y_{i_p-1})\prod_{j=1}^{p}H(i_j\to i_{j+1}-1),
$$
where we have introduced the notation
$\displaystyle \tilde{H}(1\to k):=\frac{H(1\to k)}{\prod_{i=1}^k H_1(x_i;y_{i-1})}$ and omitted the condition $i_j\in[\![1,k]\!]$ for all $j\in[\![1,p]\!]$
in order to shorten the expression. Note that this reasoning is exactly
the same as the one we carried out in \cite{BEL26} when we had $k=2$.
\\

We can now use Equation \eqref{eq:ak_prodxi*prodyi} to decompose $A_p$ in
terms of the $A_1$'s. We thus find:
\begin{equation}\label{eq:htilde}
\dpart{\tilde{H}(1\to k)}{t}=\frac{1}{\displaystyle\prod_{i=1}^k H_1(x_i;y_{i-1})}\sum_{p=2}^k\sum_{\underset{1\leqslant a,b\leqslant p}{i_1<\dots<i_p}}\frac{A_1(x_{i_a},y_{i_b-1})}{\displaystyle\prod_{j\neq a}(x_{i_j}-x_{i_a})\prod_{j\neq b}(y_{i_j-1}-y_{i_b-1})}\prod_{j=1}^{p}H(i_j\to i_{j+1}-1),
\end{equation}
where the right-hand side can also be rewritten as, by exchanging $i_a,i_b$ into $\alpha,\beta+1$ and defining the set $I:=\{i_1,\dots,i_p\}$
\begin{equation}\label{eq:htilde2}
\frac{1}{\displaystyle\prod_{i=1}^k H_1(x_i;y_{i-1})}\sum_{1\leqslant\alpha,\beta\leqslant k}\sum_{\underset{\vert I\vert\geqslant2, \alpha,\beta+1\in I}{I\subset[\![1,k]\!]}}\frac{A_1(x_{i_a},y_{i_b-1})}{\displaystyle\prod_{j\neq a}(x_{i_j}-x_{i_a})\prod_{j\neq b}(y_{i_j-1}-y_{i_b-1})}\prod_{j=1}^{p}H(i_j\to i_{j+1}-1),
\end{equation}
where we recall the cyclic convention $\beta+1\equiv 1$ if $\beta=k$.
\\

As stated at the end of the previous section, we assume that Proposition \ref{prop:recu} holds up to 
rank $k-1$, allowing us to decompose $H_j, j\in[\![1,k-1]\!]$ as a linear combination of
the $\prod_{i=1}^jH_1(x_i;y_{\sigma(i)})$, where $\sigma\in\mathfrak{P}_j$.
We thus have in particular the formulas, for $a,b\in[\![1,k]\!], b\neq a-1$
(since $H(a\to a-1)=H(1\to k)$ as seen previously, so the induction
hypothesis does not apply there):
\begin{equation}\label{eq:h_atob}
H(a\to b)=\sum_{\sigma\in\mathfrak{P}_{[\![a,b]\!]}}\!\!C_{\sigma}(a\to b)\,\prod_{i=a}^bH_1(x_{i};y_{\sigma(i)}),
\end{equation}
if $a\leqslant b$ and, otherwise
$$
H(a\to b)=\sum_{\sigma\in\mathfrak{P}_{[\![a,k]\!]\cup[\![1,b]\!]}}\hspace{-0.5cm}C_{\sigma}(a\to b)\,H_1(x_{a};y_{\sigma(a)})\cdots H_1(x_k;y_{\sigma(k)})\,H_1(x_1;y_{\sigma(1)})\cdots H_1(x_b;y_{\sigma(b)}),
$$
which we will hereafter summarize as the first formula, to simplify
notation.
\\

We now wish to plug this formula into Equation \eqref{eq:htilde}, more
precisely into the last product. Consider integers
$1\leqslant i_1<\dots<i_p\leqslant k$. Note that we have already isolated
the term for $p=1$, so there are at most $k-1$ terms with indices between
$i_j$ and $i_{j+1}-1$, which allows us to apply the induction hypothesis to
them. We thus find:
$$
\prod_{j=1}^{p}H(i_j\to i_{j+1}-1)=\prod_{j=1}^p\left(\sum_{\sigma_j\in\mathfrak{P}_{[\![i_j,i_{j+1}-1]\!]}}\hspace{-0.4cm}C_{\sigma_j}(i_j\to i_{j+1}-1)\,\prod_{i=i_j}^{i_{j+1}-1} H_1(x_i;y_{\sigma_j(i)})\right).
$$
We then use Lemma \ref{lem:permutation} to build, from our permutations
$\sigma_1,\dots,\sigma_p$ appearing in the product, a large plane
permutation $\sigma$ as a concatenation of the latter:
$$
\prod_{j=1}^{p}H(i_j\to i_{j+1}-1)=\hspace{-0.2cm}\sum_{\sigma\in\mathfrak{P}_{I,k}}\hspace{-0.2cm}C_{\sigma}(i_1\to i_2-1)\,\cdots \,C_{\sigma}(i_p\to i_1-1)\prod_{i=1}^k H_1(x_i;y_{\sigma(i)}),
$$
with $I=\{i_1,\dots,i_p\}$, and the definition of $\mathfrak{P}_{I,k}$ is
the one following Lemma \ref{lem:permutation}. Also recall that the
$C_{\sigma}$'s appearing here refer to the previous $C_{\sigma_i}$'s, since
the $\sigma_i$'s are the restrictions of $\sigma$ to the intervals
considered, see the definition of $C$ in Section \ref{sec:intro_csigma}.
\\

Now the main idea is to exchange the sum over $I=\{i_1,\dots,i_p\}\subset[\![1,k]\!]$ containing $\alpha,\beta+1$ (which are the latter $i_a,i_b$) and the one over $\sigma\in\mathfrak{P}_{I,k}$ from Equation \eqref{eq:htilde}. And 
this is possible using Proposition \ref{prop:permutation} which gives the equation $\bigcup_{\vert I\vert\geqslant 2}\mathfrak{P}_{I,k}=\mathfrak{P}_k\backslash\{\gamma_k\}$.
Hence it only remains to plug everything into Equation \eqref{eq:htilde} to
finally obtain, after a small rearrangement:
 \begin{equation*}
    \begin{split}
        \dpart{\tilde{H}(1\to k)}{t}=\displaystyle\sum_{\underset{\sigma\neq\gamma_k}{\sigma\in\mathfrak{P}_k}}\displaystyle\frac{\prod_{i=1}^kH_1(x_i;y_{\sigma(i)})}{\prod_{i=1}^k H_1(x_i;y_{i-1})}&\displaystyle\sum_{1\leqslant\alpha,\beta\leqslant k}A_1(x_{\alpha};y_{\beta})\\
        &\times \hspace{-0.4cm}\sum_{\underset{\vert I\vert\geqslant2;\, \alpha,\beta+1\in I}{I\subseteq[\![1,k]\!]/}}\hspace{-0.4cm}\frac{C_{\sigma}(I[1]\to I[2]-1)\cdots C_{\sigma}(I[\,\,\vert I\vert\,]\to I[1]-1)}{\displaystyle\prod_{\underset{i\neq\alpha}{i\in I}}(x_{i}-x_{\alpha})\prod_{\underset{j\neq\beta+1}{j\in I}}(y_{j-1}-y_{\beta})},
    \end{split}
\end{equation*}
where we observe that at least one factor $C_{\sigma}$ vanishes whenever 
$\sigma\not\in\mathfrak{P}_{I,k}$ (recall again that $C_{\sigma}(a\to b)=0$ if $[\![a,b]\!]$ is not stable under $\sigma$).\\

Finally, we may rewrite this equation by making explicit the notation
$\tilde{H}$ as well as $C_{\sigma}^{\alpha,\beta}$ introduced earlier:
$$
\displaystyle\dpart{}{t}\left(\frac{H(1\to k)}{\prod_{i=1}^k H_1(x_i;y_{i-1})}\right)=\displaystyle\sum_{\underset{\sigma\neq\gamma_k}{\sigma\in\mathfrak{P}_k}}\displaystyle\frac{\prod_{i=1}^kH_1(x_i;y_{\sigma(i)})}{\prod_{i=1}^k H_1(x_i;y_{i-1})}\displaystyle\sum_{1\leqslant\alpha,\beta\leqslant k}A_1(x_{\alpha};y_{\beta})\cdot C_{\sigma}^{\alpha,\beta}(1\to k),
$$
which is exactly the statement of Lemma \ref{lem:introcalphabeta} we wanted to show.

\subsubsection{Relation between $C_{\sigma}^{\alpha,\beta}$ and $C_{\sigma}$}

To compute these coefficients, let us begin by recalling Remark
\ref{rem:gammak}, which already allows us to compute
$C_{\gamma_k}^{\alpha,\beta}$, which is thus equal to $0$, as desired. We
can now consider $\sigma\in\mathfrak{P}_k\backslash\{\gamma_k\}$.
\\

As stated in Proposition \ref{prop:recu}, we will see two situations
appearing in our computation of the $C_{\sigma}^{\alpha,\beta}$'s depending on the values of
$\beta$. Once this is done, we can move on to the end of the proof and
conclude. So consider two integers $\alpha,\beta\in[\![1,k]\!]$. Let us
first note the equality
\begin{equation}\label{eq:sym_calphabeta}
C_{\sigma}^{\alpha,\beta}(x_1,\dots,x_k;y_1,\dots,y_k)=C_{\gamma_k^{-1}\circ\sigma\circ\gamma_k}^{\alpha+1,\beta+1}(x_2,\dots,x_k,x_1;y_2,\dots,y_k,y_1)
\end{equation}
which is proven simply by going back to the definition of
$C_{\sigma}^{\alpha,\beta}(1\to k)$ in Equation \eqref{eq:calphabeta} and then
performing the change of variable
$(I:=\{i_1,\dots,i_p\})\longmapsto \{i_1+1,\dots,i_p+1\}$ in the sum. Thanks
to this property, the rotational invariance of $C_{\sigma}$ (see the third
equation of \eqref{eq:recu_csigma}) and that we have
$\sigma\neq \gamma_k$, we can reduce to fixing $\alpha$ and choosing it such
that $\sigma(\alpha)\neq \alpha-1$. For clarity of notation, we will reduce
to the case $\alpha:=1$, which then verify $\sigma(1)\neq k$, the proof being the same
for other values of $\alpha$.
\\
Let us conclude this paragraph with a brief reminder of the definition of $C_{\sigma}^{1,\beta}(1\to n)$, defined for all $\sigma\in\mathfrak{P}_k\backslash\{\gamma_k\}$,
$\beta\in[\![1,k]\!]$ and $n\geqslant \beta$:
$$
C_{\sigma}^{1,\beta}(1\to n):=\hspace{-0.3cm}\sum_{\underset{\vert I\vert\geqslant2;\, 1,\beta+1\in I}{I\subseteq[\![1,n]\!]/}}\hspace{-0.3cm}\frac{C_{\sigma}(1\to I[2]-1)\cdots C_{\sigma}(\dots\to \beta)\,C_{\sigma}(\beta+1\to \dots) \cdots C_{\sigma}(I[\,\vert I\vert\,]\to n)}{\displaystyle\prod_{\underset{i\neq1}{i\in I}}(x_{i}-x_{1})\prod_{\underset{j\neq\beta+1}{j\in I}}(y_{j-1}-y_{\beta})}.
$$

\paragraph{Case $\beta\neq k$:}
For this case, we will highlight a specific element belonging to the
subsets $I$ over which we sum in order to isolate the last $C_{\sigma}$ of the product: their maximum, which is $I[\vert I\vert]$ by our previous conventions. Since we consider subsets
containing $1$ and $\beta+1\neq1$, we directly have
$i_0:=\max(I)-1=I[\,\vert I\vert\,]-1\geqslant \beta$. We can thus write $I=:J\cup\{i_0+1\}$ where $\vert J\vert\geqslant2$ if and only if $i_0\neq \beta$. We will hence isolate this term and also put forward all the factors where $i_0+1=\max(I)$ appears :
 \begin{equation*}
    \begin{split}
        &C_{\sigma}^{1,\beta}(1\to k)=\,\frac{C_{\sigma}(\beta+1\to k)}{(x_{\beta+1}-x_1)(y_k-y_{\beta})}\sum_{\underset{\vert J\vert\geqslant1;\, 1\in J}{J\subseteq[\![1,\beta]\!]/}}\frac{C_{\sigma}(1\to J[2]-1) \cdots \,C_{\sigma}(J[\,\vert J\vert\,]\to \beta)}{\displaystyle\prod_{\underset{i\neq1}{i\in J}}(x_{i}-x_{1})\prod_{\underset{j\neq1}{j\in J}}(y_{j-1}-y_{\beta})}\\
        &+\sum_{i_0=\beta+1}^{k-1}\frac{C_{\sigma}(i_0+1\to k)}{(x_{i_0+1}-x_1)(y_k-y_{\beta})(y_{i_0}-y_{\beta})}\sum_{\underset{\vert J\vert\geqslant2;\, 1,\beta+1\in J}{J\subseteq[\![1,i_0]\!]/}}\frac{C_{\sigma}(1\to J[2]-1) \cdots \,C_{\sigma}(J[\,\vert J\vert\,]\to i_0)}{\displaystyle\prod_{\underset{i\neq1}{i\in J}}(x_{i}-x_{1})\prod_{\underset{j\neq\beta+1,1}{j\in J}}(y_{j-1}-y_{\beta})},
    \end{split}
\end{equation*}
where the first term corresponds to the case $i_0=\beta$ and the second to the case $i_0\neq \beta$.
We can now bring the term $y_{i_0}-y_{\beta}$ back into the product, since
it corresponds to the case $j=1$ (the elements of $J\subseteq[\![1,i_0]\!]$ are now considered
modulo $i_0$ rather than $k$). We have thus the relation:
 \begin{equation*}
    \begin{split}
        C_{\sigma}^{1,\beta}(1\to k)=&\,\frac{C_{\sigma}(\beta+1\to k)}{(x_{\beta+1}-x_1)(y_k-y_{\beta})}\sum_{\underset{\vert J\vert\geqslant1;\, 1\in J}{J\subseteq[\![1,\beta]\!]/}}\frac{C_{\sigma}(1\to J[2]-1) \cdots \,C_{\sigma}(J[\,\vert J\vert\,]\to \beta)}{\displaystyle\prod_{\underset{i\neq1}{i\in J}}(x_{i}-x_{1})\prod_{\underset{j\neq1}{j\in J}}(y_{j-1}-y_{\beta})}\\
        &+\sum_{i_0=\beta+1}^{k-1}\frac{C_{\sigma}(i_0+1\to k)}{(x_{i_0+1}-x_1)(y_k-y_{\beta})}\sum_{\underset{\vert J\vert\geqslant2;\, 1,\beta+1\in J}{J\subseteq[\![1,i_0]\!]/}}\frac{C_{\sigma}(1\to J[2]-1) \cdots \,C_{\sigma}(J[\,\vert J\vert\,]\to i_0)}{\displaystyle\prod_{\underset{i\neq1}{i\in J}}(x_{i}-x_{1})\prod_{\underset{j\neq\beta+1}{j\in J}}(y_{j-1}-y_{\beta})}.
    \end{split}
\end{equation*}
By highlighting the term $\vert J\vert=1$ (which is $J=\{1\}$), the first sum over $J$ is finally equal to $C_{\sigma}(1\to \beta)+C_{\sigma}^{1,\beta}(1\to\beta)$ while the last sum over $J$ equals
$C_{\sigma}^{1,\beta}(1\to i_0)$, leading to 
$$
C_{\sigma}^{1,\beta}(1\to k)=\frac{(C_{\sigma}(1\to \beta)+C_{\sigma}^{1,\beta}(1\to\beta))\,C_{\sigma}(\beta+1\to k)}{(x_{\beta+1}-x_1)(y_k-y_{\beta})}+\sum_{i_0=\beta+1}^{k-1}\frac{C_{\sigma}^{1,\beta}(1\to i_0)\,C_{\sigma}(i_0+1\to k)}{(x_{i_0+1}-x_1)(y_k-y_{\beta})}.
$$
One can then use the induction hypothesis, since $\beta\leqslant i_0<k$, in order to compute the coefficient $C_{\sigma}^{1,\beta}(1\to i_0)$. 
First, consider the case $\beta\neq\sigma(1)$, hence $C_{\sigma}(1\to \beta)+C_{\sigma}^{1,\beta}(1\to\beta)=0$ and $C_{\sigma}^{1,\beta}(1\to i_0)=C_{\sigma}(1\to i_0)$ for all  $i_0\in[\![\beta+1,k-1]\!]$. Then everything vanishes and this case is done.\\

Let us now take $\beta=\sigma(1)$, which leads to $C_{\sigma}^{1,\beta}(1\to\beta)=0$ and $C_{\sigma}^{1,\beta}(1\to i_0)=C_{\sigma}(1\to i_0)$ for all $i_0\in[\![\beta+1,k-1]\!]$. We can write
$$
C_{\sigma}^{1,\beta}(1\to k)=\frac{C_{\sigma}(1\to \beta)\,C_{\sigma}(\beta+1\to k)}{(x_{\beta+1}-x_1)(y_k-y_{\beta})}+\sum_{i_0=\beta+1}^{k-1}\frac{C_{\sigma}(1\to i_0)\,C_{\sigma}(i_0+1\to k)}{(x_{i_0+1}-x_1)(y_k-y_{\beta})},
$$
since the first term is nothing but the case $i_0=\beta$ in the sum. However $[\![1,i_0]\!]$ is not stable under $\sigma$ since $\sigma(1)=\beta\not\in[\![1,i_0]\!]$, and thus
$C_{\sigma}(1\to i_0)=0$ in this case. We can thus extend the sum as:
$$
C_{\sigma}^{1,\beta}(1\to k)=\sum_{i_0=1}^{k-1}\frac{C_{\sigma}(1\to i_0)\,C_{\sigma}(i_0+1\to k)}{(x_{i_0+1}-x_1)(y_k-y_{\beta})},
$$
which is exactly the formula defining $C_{\sigma}(1\to k)$ for
$\sigma(1)\neq k$ seen in Equation \eqref{eq:recu_csigma}!

\paragraph{Case $\beta=k$.}

In this case, we will use a small trick. Consider a complex polynomial $P$ such that
$P(X)=\prod_{j=1}^p(z_j-X)$ where the $z_i$'s are distinct, $p\geqslant2$.
Then we have the relation:
$$
\sum_{i=1}^p\frac{1}{P'(z_i)}=-\sum_{i=1}^p\frac{1}{\prod_{j\neq i}(z_j-z_i)}=0,
$$
as a sum of residues of the meromorphic function $\frac{1}{P}$ defined on
the sphere $\C\cup\{\infty\}$. This property naturally generalizes to
formal variables $z_i$ by considering a partial fraction decomposition in
the appropriate ring.
\\

Thus, considering $P=P_I:=\prod_{j\in I}(y_{j-1}-X)$ for any subset
$I\subseteq[\![1,k]\!]$, the definition of the $C_{\sigma}^{\alpha,\beta}$'s
together with the property we have just shown translate into the equation
$$
\sum_{\beta=1}^kC_{\sigma}^{1,\beta}(1\to k)=\hspace{-0.3cm}\sum_{\underset{\vert I\vert\geqslant2;\, 1\in I}{I\subseteq[\![1,n]\!]/}}\frac{C_{\sigma}(1\to I[2]-1)\cdots C_{\sigma}(I[\,\vert I\vert\,]\to k)}{\displaystyle\prod_{\underset{i\neq1}{i\in I}}(x_{i}-x_{1})}\sum_{\beta/\beta+1\in I}\frac{1}{P_I'(y_{\beta})}
=0.
$$
We can now use what we have already shown for the case $\beta\neq k$ to
finally deduce that
$$
0=\sum_{\beta=1}^kC_{\sigma}^{1,\beta}(1\to k)=\underbrace{C_{\sigma}^{1,\sigma(1)}(1\to k)}_{=C_{\sigma}(1\to k)}+C_{\sigma}^{1,k}(1\to k)+\underbrace{\sum_{\beta\neq \sigma(1),k}C_{\sigma}^{1,\beta}(1\to k)}_{=0},
$$
hence $C_{\sigma}^{1,k}(1\to k)=-C_{\sigma}(1\to k)$. This concludes
our computation of the $C_{\sigma}^{\alpha,\beta}$s and the proof of the first statement of Proposition \ref{prop:recu}!

\subsubsection{End of the proof}\label{sec:conclusion_hk}

Now that we have computed the coefficients $C_{\sigma}^{\alpha,\beta}$, it
only remains to plug these data into Lemma \ref{lem:introcalphabeta} to deduce the second statement of Proposition \ref{prop:recu}. We
then find:
\begin{equation*}
    \begin{split}
\displaystyle\dpart{}{t}\left(\frac{H(1\to k)}{\prod_{i=1}^k H_1(x_i;y_{i-1})}\right)=
\displaystyle\sum_{\underset{\sigma\neq\gamma_k}{\sigma\in\mathfrak{P}_k}}C_{\sigma}(1\to k)\,\,
\displaystyle&\frac{\prod_{i=1}^kH_1(x_i;y_{\sigma(i)})}{\prod_{i=1}^k H_1(x_i;y_{i-1})}\\
&\times\displaystyle\left(\sum_{i=1}^kA_1(x_i;y_{\sigma(i)})-\sum_{i=1}^kA_1(x_i;y_{i-1})\right).
    \end{split}
\end{equation*}
This formula can be simplified using Equation \eqref{eq:dH1/dt}, which lets
us recognize a logarithmic derivative and thus write:
$$
\displaystyle\frac{\prod_{i=1}^kH_1(x_i;y_{\sigma(i)})}{\prod_{i=1}^k H_1(x_i;y_{i-1})}\cdot\left(\sum_{i=1}^kA_1(x_i;y_{\sigma(i)})-\sum_{i=1}^kA_1(x_i;y_{i-1})\right)=\displaystyle\dpart{}{t}\left(\displaystyle\frac{\prod_{i=1}^kH_1(x_i;y_{\sigma(i)})}{\prod_{i=1}^k H_1(x_i;y_{i-1})}\right).
$$
But the $C_{\sigma}$'s do not depend on $t$ by construction, so we can
integrate the previous equation to find
$$
\left(\frac{H(1\to k)}{\prod_{i=1}^k H_1(x_i;y_{i-1})}\right)=K+\displaystyle\sum_{\underset{\sigma\neq\gamma_k}{\sigma\in\mathfrak{P}_k}}C_{\sigma}(1\to k)\displaystyle\frac{\prod_{i=1}^kH_1(x_i;y_{\sigma(i)})}{\prod_{i=1}^k H_1(x_i;y_{i-1})},
$$
where $K$ is a constant independent of $t$. To compute it, we now refer to
the definition of the $H_k$'s, seen in Equation \eqref{eq:Hkdef}, which
allows us to notice that $[t^0]H_k=\delta_{k,1}$ and thus compute $K$. We
find:
$$
K=-\displaystyle\sum_{\underset{\sigma\neq\gamma_k}{\sigma\in\mathfrak{P}_k}}C_{\sigma}(1\to k)=C_{\gamma_k}(1\to k),
$$
since the definition of $C_{\gamma_k}$ is exactly such that
the sum of the $C_{\sigma}, \sigma\in\mathfrak{P}_k$ vanishes, cf Equation \eqref{eq:recu_csigma}. It only
remains to multiply by $\prod_{i=1}^k H_1(x_i;y_{i-1})$ to finally obtain:
\begin{align*}
H(1\to k)&=C_{\gamma_k}(1\to k)\,\prod_{i=1}^k H_1(x_i;y_{i-1})+\displaystyle\sum_{\underset{\sigma\neq\gamma_k}{\sigma\in\mathfrak{P}_k}}C_{\sigma}(1\to k)\displaystyle\prod_{i=1}^kH_1(x_i;y_{\sigma(i)})\\
&=\sum_{\sigma\in\mathfrak{P}_k}C_{\sigma}(1\to k)\prod_{i=1}^kH_1(x_i;y_{\sigma(i)}).
\end{align*}
This concludes the proof of Proposition \ref{prop:recu} and thus that of
Corollary \ref{cor:hk}.

\section{A formula for strongly connected hypermaps}\label{sec:bk}

In this section we prove Theorem \ref{th:bk}, which
concerns the generating function of strongly connected alternating
hypermaps — that is, alternating hypermaps with
no bridge. The proof proceeds in three steps: we first express the
generating function of general hypermaps in terms of that of strongly
connected ones; we then invert this relation to express the latter in
terms of the former; finally, we manipulate the resulting relation to
obtain the concise formula stated in the theorem.
\\

Before we begin, we introduce a few notations.

\subsection{A few introductory notations}

\begin{figure}[h!]
    \centering
    \includegraphics[width=0.45\textwidth]{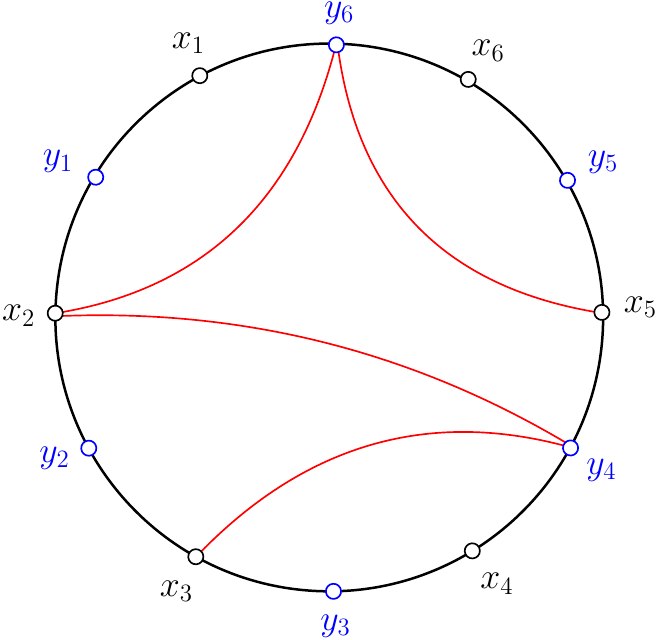}
    \caption{A plane graph $G_I$ associated with the set $I:=\{(2,4), (2,6), (3,4), (5,6)\}\subseteq E_6$.\\ Constrained edges from $E_k^c:=[\![1,k]\!]^2\backslash E_k$ are shown black and edges from $I$ are red.}
    \label{fig:graphe_GI}
\end{figure}

Let $k\geqslant2$ be an integer and consider the sets $E_k:=\{(i,j)\in[\![1,k]\!]^2, j\neq i, i-1\}$ and its complement $E_k^c:=\{(i,j)\in[\![1,k]\!]^2, j\in\{i, i-1\}\}$. As in the definition of a
plane permutation given in Section \ref{sec:planaire}, we associate to
each $I\subseteq E_k$ a bipartite graph $G_I$, whose
vertices are labeled cyclically $x_1,y_1,\dots, x_k,y_k$, with edges of
the form $x_i\sim y_j$ for $(i,j)\in I$.
\\
To every such graph $G_I$ we also add all the edges
$y_{i-1}\sim x_i \sim y_{i}$ coming from $E_k^c$.
As shown on Figure \ref{fig:graphe_GI}, the graph $G_I$ is naturally drawn in 
the plane by placing the vertices along a circle and adding chords corresponding to the edges in $I$.
 Note that, unlike a graph associated with a plane
permutation, these graphs are simple and the degree of a vertex can be arbitrary. Also, these graphs 
are allowed to be non-planar but we will be interested in the case where they are.

When $G_I$ is plane, it can be viewed as a plane map, and we may
therefore consider its faces (only the internal ones, since the outer
face is of no interest to us here). By construction, the vertices of
each face can be read cyclically as
$x_{i_1} \sim y_{j_1} \sim\dots \sim x_{i_p} \sim y_{j_p}\sim x_{i_1}$.
To every face $F$ of this form we assign the quantity
$$
H(F):=H_p(x_{i_1},\dots,x_{i_p};y_{j_1},\dots,y_{j_p}),
$$
and similarly for $B(F)$. Note that the choice of starting point and
direction in reading the vertices of the face has no effect on this
definition, since our generating series are invariant under rotation
and flip, which generate the dihedral group $D_k$.
\\
Finally, we write $h$ for the function that assigns to every edge
$e:= x_i\sim y_j$ of $[\![1,k]\!]^2$ the series $h(e):=H_1(x_i;y_j)$.
\\

We now have all the notations needed to prove the propositions that
follow.

\subsection{From blobs to general hypermaps}

As said above, every strongly connected $k$-alternating
hypermap is one that contains no bridge. The purpose of this section 
is to find a relation between $H(1\to k)$ and the $B_i$'s, the idea being to cut all 
the bridges of an arbitrary $k$-alternating hypermap in order to decompose 
it into its strongly connected components counted by the $B_i$'s. 
Such decomposition can be represented as a plane graph $G_I$, the type considered in the 
previous subsection, where the edges of $G_I$ 
encode the position of the bridges of the hypermap. 
\\

\begin{figure}[h!]
    \centering
    \includegraphics[width=1\textwidth]{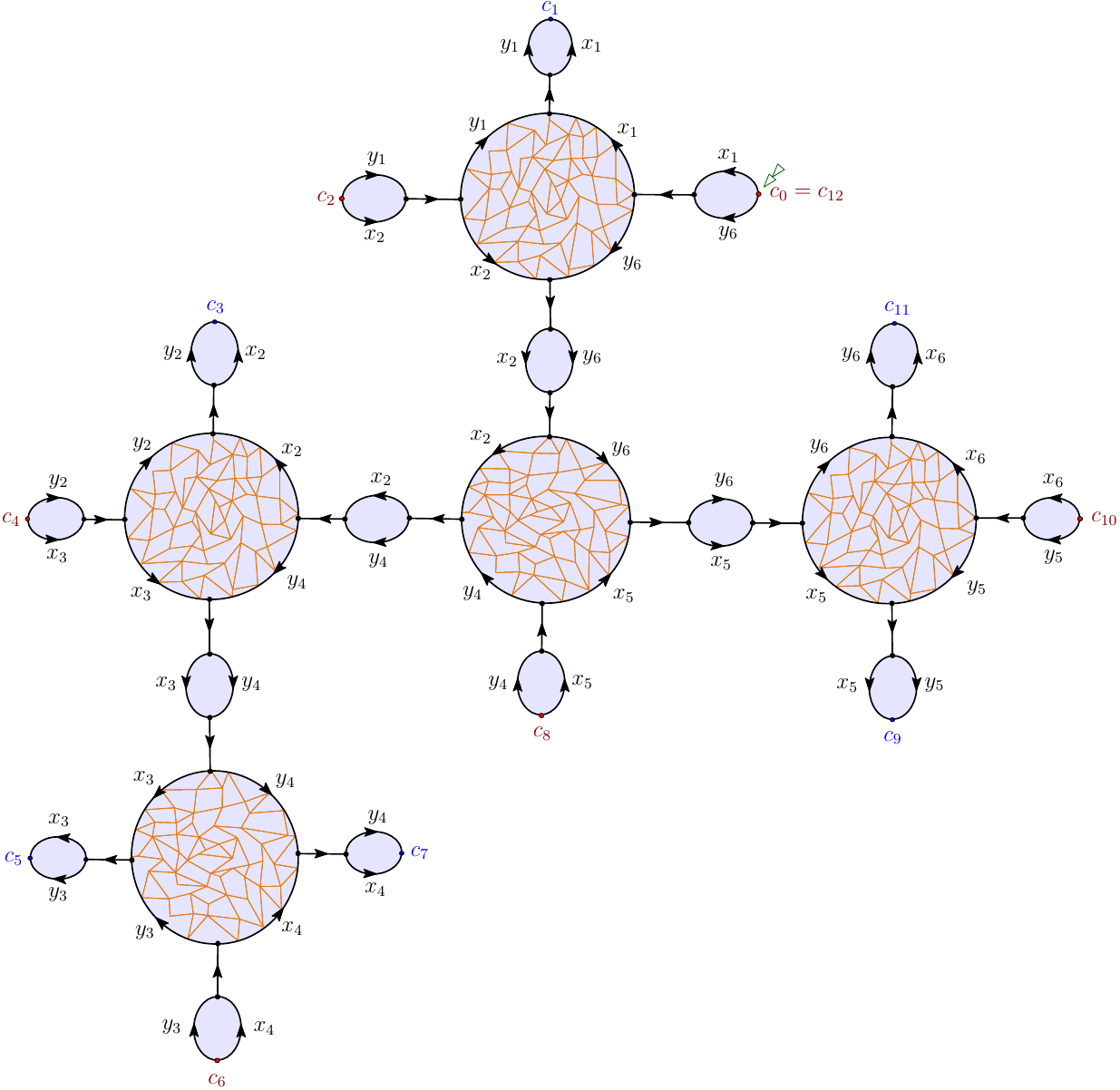}
    \caption{General structure of a $6$-alternating hypermap associated with the plane graph from Figure \ref{fig:graphe_GI}, we recall Figures \ref{fig:schema_sans_rien} and Figure \ref{fig:schema_sans_rien_stronglyconnected} for the 
    sketch conventions. Each black or red edge $x_i\sim y_j$ from $G_I$ corresponds to a $1$-alternating hypermap counted by $H_1(x_i,y_j)$, and each face of $G_I$ corresponds to a blob with its associated weights displayed on the boundary.}
    \label{fig:blob_decompo}
\end{figure}

More precisely, to each $k$-alternating hypermap, we associate the set $I$ of all pairs $(i,j)\in E_k$
such that the hypermap contains a bridge of weight $x_i^{-1}y_j^{-1}$. It turns out that the corresponding 
graph $G_I$ is plane since the hypermap is, and that the generating function $z_I$ of 
$k$-alternating hypermaps associated to a prescribed plane graph $G_I$ reads
\begin{equation}\label{eq:zI}
z_I=\left(\prod_{F\,\text{face of $G_I$}}\hspace{-0.3cm}B(F)\right)\left(\prod_{e\,\text{edge of $G_I$}}\hspace{-0.3cm}h(e)\right),
\end{equation}
where $B(F)$ and $h(e)$ are defined in the previous subsection.
The reason why such a formula holds is that the hypermaps contributing to $z_I$ have a common 
general structure, illustrated in Figure \ref{fig:blob_decompo}, which we now explain in more detail. 
\\

Let us fix a plane graph $G_I, I\subset E_k$ and consider a pair $(i,j)\in I$.
Any hypermap contributing to $z_I$ contains by definition a bridge of weight
$x_i^{-1}y_j^{-1}$. Note that several
bridges of this weight may occur. However, any two such consecutive bridges delimit a
strongly connected $1$-alternating hypermap counted by $B_1(x_i;y_j)$. Hence the first and last bridges delimit
a sequence of such $1$-alternating blobs, i.e., a $1$-alternating hypermap of weight
$H_1(x_i;y_j)$ (recall Formula \eqref{eq:h1b1} here). Note that, if there is 
a single bridge of weight $x_i^{-1}y_j^{-1}$ then the sequence is empty, which 
corresponds to the constant term $1$ appearing in $H_1$. We accordingly assign to each edge
$e:=x_i\sim y_j$ of our graph the weight $H_1(x_i;y_j)=:h(e)$. 
\\

Note that the cases $j=i$ and $j=i-1$ are slightly different,
since here a single bridge — rather than two — already implies the existence
of a blob of weight $B_1(x_i,y_j)$. Still, we may instead treat these edges as if a
 weight $H_1(x_i;y_j)$ is always
present, already assuming these edges belong to our graph; the term
$1$ then accounts for the case where there is no bridge of weight
$x_i^{-1} y_j^{-1}$. This is precisely what justifies
the particular definition of the set $E_k$.
\\

Now that we understand where the factor $\prod_{e\,\text{edge of $G_I$}}h(e)$
 comes from in Equation \eqref{eq:zI}, it remains only to explain why each 
 face $F$ of our graph yields a factor $B(F)$.
Since each bridge is associated to an edge of our graph, 
a face of $G_I$ corresponds to a bridge-free region — that is, by definition, 
a strongly connected hypermap —  of each hypermap contributing to $z_I$.
More precisely, a face $F$ with incident vertices 
$x_{i_1} \sim y_{j_1} \sim\dots \sim x_{i_p} \sim y_{j_p}\sim x_{i_1}$, corresponds
to a strongly connected $p$-alternating hypermap counted by the generating function 
of blobs $B_p(x_{i_1},\dots,x_{i_p};y_{j_1},\dots,y_{j_p})=:B(F)$.
\\

This establishes the wanted expression \eqref{eq:zI} for $z_I$; note that it is not 
necessary to add any monomial weight in the $x_i^{-1}$ and 
$y_j^{-1}$ as the weights of the bridges are 
already properly taken into account thanks to the extra 
weight factor $\prod_{i=1}^p x_i^{-1}y_i^{-1}$ of each $p$-alternating hypermap as defined in 
Equation \eqref{eq:weight}.

Now observing that $H(1\to k)$ is obtained by summing all the $z_I$ such that $G_I$ is plane, we get the following proposition. 
\begin{prop}\label{prop:hkbk}
Using the notations of this section and the previous one, we can relate
the generating series $H(1\to k)$ to that of strongly connected
alternating hypermaps by:
\begin{equation}\label{eq:hkbi}
H(1\to k)=\sum_{\underset{G_I\,\,\text{plane}}{I\subseteq E_k /}}\left(\prod_{F\,\text{face of $G_I$}}\hspace{-0.3cm}B(F)\right)\left(\prod_{e\,\text{edge of $G_I$}}\hspace{-0.3cm}h(e)\right).
\end{equation}
Recalling that every graph $G_I, I\subseteq E_k$ already contains the
edges $y_{i-1}\sim x_i\sim y_{i}$, this can equivalently be written as:
\begin{equation}\label{eq:hkbk}
\frac{H(1\to k)}{\prod_{i=1}^kH_1(x_i;y_i)H_1(x_i;y_{i-1})}=\sum_{\underset{G_I\,\,\text{plane}}{I\subseteq E_k}}\left(\prod_{F\,\text{face of $G_I$}}\hspace{-0.3cm}B(F)\right)\left(\prod_{e\in I}\,h(e)\right).
\end{equation}
\end{prop}

The statement of this proposition can be generalized as follows, which
will be useful later. Let $P\subseteq E_k$ be a fixed set, and consider the generating series 
of $k$-alternating hypermaps whose associated 
graph contains (at
least) the edges in $P$.
On the one hand, this quantity is obtained by 
summing the $z_I$ over all $I\subseteq E_k$ such that $P\subseteq I$ and $G_I$ is a plane graph.
On the other hand, using a variant of the reasoning leading to Equation \eqref{eq:zI}
 where we only decompose along the bridges corresponding to the elements of $P$,
we find that our quantity is nothing but the product of the generating series $H(F)$ 
over all the faces of the graph $G_{P}$, divided by the product of $h(e)$ over all the edges $e$ of $G_P$, which are $P\sqcup E_k^c$.
Indeed, the edges $e$ lying on the boundary of each face of $G_P$ —
belonging to the set $P\sqcup E_k^c$— can no longer take part in the
subgraphs, so a weight $h(e)$ must be removed for each of them. We end up 
with the identity
\begin{equation}\label{eq:hkbk2}
\frac{\prod_{\underset{\text{of $G_P$}}{\text{$F$ face}}}H(F)}{\displaystyle\prod_{e\in P\sqcup E_k^c}h(e)}=\sum_{\underset{G_I\,\,\text{plane}}{I\colon P\subseteq I\subseteq E_k}}\left(\prod_{F\,\text{face of $G_I$}}\hspace{-0.3cm}B(F)\right)\left(\prod_{e\in I}\,h(e)\right).
\end{equation}
Formula \eqref{eq:hkbk} is therefore simply the case $P=\emptyset$ of
this formula, recalling that $G_{\emptyset}$ is just the graph made up
of the edges $x_1\sim y_1 \sim x_2 \sim \dots \sim y_k\sim x_1$.

\subsection{From general hypermaps to blobs via inclusion/exclusion}

Interestingly, the Relation \eqref{eq:hkbk2} may be inverted using the inclusion/exclusion principle, which we invoke in 
the following form.

\begin{lem}[Inclusion/Exclusion principle]\label{lem:incluexclu}
Let $E$ be a finite set and $f\colon \mathcal{P}(E)\longrightarrow \C$ a function. We also define the function $f_{-1}\colon \mathcal{P}(E)\longrightarrow \C$ by:
$$
\forall P\subseteq E, f_{-1}(P):=\sum_{I\colon  P\subseteq I\subseteq E}f(I).
$$
Then the following \textbf{inclusion/exclusion formula} holds:
\begin{equation}\label{eq:incexc}
\forall P\subseteq E, f(P)=\sum_{I\colon P\subseteq I\subseteq E}(-1)^{\vert I\vert -\vert P\vert}\,f_{-1}(I).
\end{equation}
\end{lem}

Let us recall for completeness the proof which is short and follows from a straightforward computation.

\begin{proof}
For $P\subseteq E$, we have
\begin{align*}
\sum_{I\colon P\subseteq I\subseteq E}(-1)^{\vert I\vert -\vert P\vert}\,f_{-1}(I)&=\sum_{(K,I)\colon P\subseteq I\subseteq K\subseteq E}(-1)^{\vert I\vert -\vert P\vert}\,f(K)\\
&=\sum_{K\colon P\subseteq K\subseteq E}f(K)\sum_{I\colon \,P\subseteq I\subseteq K}(-1)^{\vert I\vert -\vert P\vert}\\
&\overset{"i=\vert I\backslash P\vert"}{=}\!\!\!\sum_{K\colon P\subseteq K\subseteq E}f(K)\sum_{i=0}^{\vert K\vert-\vert P\vert}(-1)^{i}\,\binom{\vert K\vert -\vert P\vert}{i}\\
&=\sum_{K\colon P\subseteq K\subseteq E}f(K)\,\,\mathbf{1}_{\vert P\vert-\vert K\vert,0}\\
&=f(P).
\end{align*}
\end{proof}

To apply our newly acquired lemma, we take $E=E_k$ and observe that Equation \eqref{eq:hkbk2} may be rewritten in the form:
$$
\forall P\subseteq E_k, f_{-1}(P,h):=\sum_{I \colon P\subseteq I\subseteq E_k}f(I,h),
$$
with
\begin{equation}\label{eq:f1}
f(I,h):=\mathbf{1}_{G_I\,\,\text{plane}}\hspace{-0.3cm}\prod_{F\,\text{face of $G_I$}}\hspace{-0.5cm}B(F)\,\,\,\prod_{e\in I}h(e),
\end{equation}
and
\begin{equation}\label{eq:f-1}
f_{-1}(P,h)=\sum_{I \colon P\subseteq I\subseteq E_k}f(I,h)=\mathbf{1}_{G_P\,\,\text{plane}}
\displaystyle\,\frac{\prod_{\underset{\text{of $G_P$}}{\text{$F$ face}}}H(F)}{\displaystyle\prod_{e\in P\sqcup E_k^c}h(e)}.
\end{equation}
Note that we added the extra parameter $h$ in those equation that will be useful later.
By the Inclusion/Exclusion Lemma \ref{lem:incluexclu}, we may conversely express $f$ in terms of $f_{-1}$ as:
\begin{equation}\label{eq:f12}
f(P,h)=\sum_{J\colon P\subseteq J\subseteq E_k}(-1)^{\vert J\vert -\vert P\vert}\,f_{-1}(J,h)=\sum_{J\colon P\subseteq J\subseteq E_k}(-1)^{\vert J\vert -\vert P\vert}\,\,\mathbf{1}_{G_J\,\,\text{plane}}\,
\frac{\prod_{\underset{\text{of $G_J$}}{\text{$F$ face}}}H(F)}{\displaystyle\prod_{e\in J\sqcup E_k^c}h(e)}.
\end{equation}
Taking $P=\emptyset$, we obtain the relation:
$$
B(1\to k)=f(\emptyset,h)=\sum_{\underset{G_J\,\,\text{plane}}{J\subseteq E_k}}(-1)^{\vert J\vert}\,\,
\frac{\prod_{\underset{\text{of $G_J$}}{\text{$F$ face}}}H(F)}{\displaystyle\prod_{e\in J\sqcup E_k^c}h(e)}.
$$
This relation can be rewritten, noting that there are exactly $\vert J\vert+1$ inner faces in $G_J$ since $G_J$ is plane:
\begin{equation}\label{eq:bkhk}
\boxed{-B(1\to k)=\sum_{\underset{G_J\,\,\text{plane}}{J\subseteq E_k}}
\frac{\prod_{\underset{\text{of $G_J$}}{\text{$F$ face}}}(-H(F))}{\displaystyle\prod_{e\in J\sqcup E_k^c}h(e)}.}
\end{equation}
Note the close resemblance between Equation \eqref{eq:bkhk} and Equation \eqref{eq:hkbi}: it appears that we have simply swapped $h$ for $1/h$ and $B$ for $-H$. This correspondence is precisely what we focus on in the following subsection.

\subsection{A subsidiary equation between $H_k$ and $B_k$}

In this subsection we start from Equation \eqref{eq:bkhk}, our point will be to understand how swapping $h$ to $1/h$ will enable us to find a new relation directly between $B_k$ and $H_k$ (and no more between $B_k$ and the $H_i$'s, $i\leqslant k$).
To do so, \textbf{suppose, by induction, that Theorem} \ref{th:bk} \textbf{holds up to rank} $\mathbf{k-1, k\geqslant 3}$. The base case starts at $k=2$ and is already covered by Relation \eqref{eq:b2}, proved in \cite{BEL26}, so we focus on the induction step. 
\\
If one defines $H_p(F,h)$, for all $p\leqslant k$ and face $F$ with incident vertices 
$x_{i_1} \sim y_{j_1} \sim\dots \sim x_{i_p} \sim y_{j_p}\sim x_{i_1}$, as
\begin{equation}\label{eq:hf}
H(F,h):=\sum_{\sigma\in\mathfrak{S}_p}C_{\sigma}(F)\prod_{a=1}^ph((i_a,y_{j_{\sigma(a)}}))
\end{equation}
 where we recall that the notation $C_{\sigma}(F)$ follows the same convention as $H(F)$ and $B(F)$: it simply denotes the function $C_{\sigma}$ evaluated at the variables encountered while going around the face $F$ counterclockwise, starting from an $x$-vertex. Rotational invariance of $C_{\sigma}$ means the choice of starting point plays no role. 
 And we recall that $h((i,j)):=H_1(x_i;y_j)$ such that this function $h$ is the only parameter that gathers all the physics of alternating hypermaps. Similarly, we define $B_p(F,h)$ for all $p<k$ (the case $p=k$ remains to be proved):
\begin{equation}
B(F,h):=\sum_{\sigma\in\mathfrak{S}_p}C_{\sigma}(F)\prod_{a=1}^ph((i_a,y_{j_{\sigma(a)}}))^{-1}.
\end{equation}
Note that these definitions coincide with the formulas we know for $H(F)$ and $B(F)$ since it is no more than a change of indices in the boundary variables appearing in Corollary \ref{cor:hk} and Theorem \ref{th:bk}.
\\

Hence these relations may also be viewed as
 \begin{equation}\label{eq:bfhf}
 B(F,h)=-H(F,h^{-1}),
 \end{equation}
 for all faces $F$ of degree $p<k$. Now if we denote as $\emptyset_k$ the face with incident vertices 
$x_{1} \sim y_{1} \sim\dots \sim x_{k} \sim y_{k}\sim x_{1}$, one can define $B(1\to k)=:B(\emptyset_k,h)$ thanks to Equation \eqref{eq:bkhk}, by expressing explicitly $H(F)$ in terms of $h$. Then our purpose is to show the equality $B(\emptyset_k,h)=-H(\emptyset_k,h^{-1})$. \\

To do so, let us isolate, in Equation \eqref{eq:hkbk}, the term $I=\emptyset$ (associated with the face $\emptyset_k$) in the right-hand side. This gives:
$$
\frac{H(\emptyset_k,h)}{\prod_{i=1}^kH_1(x_i;y_i)H_1(x_i;y_{i-1})}-B(\emptyset_k,h)=\sum_{\underset{G_J\,\,\text{plane}}{\emptyset\,\neq I\subseteq E_k}}\left(\prod_{F\,\text{face of $G_I$}}\hspace{-0.3cm}B(F)\right)\left(\prod_{e\in I}\,h(e)\right).
$$
which, by multiplying by $\prod_{i=1}^kH_1(x_i;y_i)H_1(x_i;y_{i-1})$, which equals $\prod_{e\in E_k^c}h(e)$, can also be written as
\begin{equation}\label{eq:bononsenfiche}
H(\emptyset_k,h)-B(\emptyset_k,h)\prod_{e\in E_k^c}h(e)=\sum_{\underset{G_I\,\,\text{plane}}{\emptyset\,\neq I\subseteq E_k}}\left(\prod_{F\,\text{face of $G_I$}}\hspace{-0.3cm}B(F,h)\right)\left(\prod_{e\in  I\sqcup E_k^c}h(e)\right).
\end{equation}
since the degree of each face appearing in the right-hand side is at most $k-1$, so $B(F)=B(F,h)$.
\\

Now let us go back to Equation \eqref{eq:bkhk} and, similarly, isolate the $I=\emptyset$ term from the right-hand side:
$$
-B(\emptyset_k,h)+\frac{H(\emptyset_k,h)}{\prod_{e\in E_k^c}h(e)}=\sum_{\underset{G_I\,\,\text{plane}}{\emptyset\,\neq I\subseteq E_k}}
\frac{\prod_{\underset{\text{of $G_I$}}{\text{$F$ face}}}\,(-H(F,h))}{\displaystyle\prod_{e\in  I\sqcup E_k^c}h(e)}.
$$
Let us rewrite this new equation with $h$ replaced by $h^{-1}$:
$$
-B(\emptyset_k,h^{-1})+H(\emptyset_k,h^{-1})\,\prod_{e\in E_k^c}h(e)=\sum_{\underset{G_I\,\,\text{plane}}{\emptyset\,\neq I\subseteq E_k}}
\prod_{\underset{\text{of $G_I$}}{\text{$F$ face}}}\,(-H(F,h^{-1}))\cdot\left(\displaystyle\prod_{e\in  I\sqcup E_k^c}h(e)\right),
$$
and, by applying the induction hypothesis, which is Equation \eqref{eq:bfhf}, on the right-hand side, we get:
$$
-B(\emptyset_k,h^{-1})+H(\emptyset_k,h^{-1})\prod_{e\in E_k^c}h(e)=\sum_{\underset{G_I\,\,\text{plane}}{I\subseteq E_k}}
\left(\prod_{\underset{\text{of $G_I$}}{\text{$F$ face}}}\,B(F,h)\right)\left(\prod_{e\in  I\sqcup E_k^c}h(e)\right).
$$
The right-hand side is thus exactly the same as the right-hand side from Equation \eqref{eq:bononsenfiche} ! So we may write the equality:
$$
H(\emptyset_k,h)-B(\emptyset_k,h)\prod_{e\in E_k^c}h(e)=-B(\emptyset_k,h^{-1})+H(\emptyset_k,h^{-1})\prod_{e\in E_k^c}h(e),
$$
or equivalently
\begin{equation}\label{eq:bh_bhmoinsun}
H(\emptyset_k,h)+B(\emptyset_k,h^{-1})=(B(\emptyset_k,h)+H(\emptyset_k,h^{-1}))\prod_{e\in E_k^c}h(e).
\end{equation}
which is a beautiful equation that directly link $B_k$ and $H_k$ and will be useful later on.

\subsection{Decomposing $B_k$ as a function of $H_1$ via an involution}

For now, we still do not have a clear idea of why $B_k$ can be related to $H_1$ 
as the same form as $H_k$ is. We only know that $B_k$ depends on the face coefficients $t_i^{\circ}, t_j^{\bullet}$ as a function of $h$ but this is not enough. This will be the goal of this subsection.

\paragraph{Rewriting Equation \eqref{eq:bkhk}.} Along the same lines as the previous subsection, 
we start again from Relation \eqref{eq:bkhk}, expanding the $H(F)$'s terms using Equation \eqref{eq:hf}, leading to
 $$
 B(1\to k)=\sum_{\underset{G_j\,\,\text{plane}}{J\subseteq E_k}}(-1)^{\vert J\vert}\frac{1}{\prod_{e\in J\sqcup E_k^c}h(e)}\prod_{\underset{\text{de } G_P}{F\,\,\text{face}}}\left(\sum_{\sigma_F\in\mathfrak{P}_F}C_{\sigma_F}(F)\prod_{\underset{\text{de } F}{e\text{ arête}}}h(e)\right),
 $$
 The first idea one may have now is to build a single global plane permutation 
 on $[\![1,k]\!]$ from the plane permutations defined on each face. But since several edges may leave the same vertex, it will lead either to multi-valued images for a single preimage, or to a failure of injectivity: neither of which makes sense here in order to obtain a permutation. 
We therefore will gather all these informations, for a fixed set $P\subseteq E_k$ and a sequence of plane permutations $(\sigma_F)_{F}$ indexed by the faces of $G_P$, in a function $g\colon [\![1,k]\!]\longrightarrow \mathcal{T}([\![1,k]\!])$ where $\mathcal{T}([\![1,k]\!])$ denotes the set of tuples of elements of $[\![1,k]\!]$. 
We define this function such that, if we set a total order over all possible faces of plane graphs $G_I, I\subseteq E_k$, then $g(i)$ will be the sequence of all images of the $\sigma_F(i)$ where $F$ browse the total order, and $0$ for all faces that does not exist in $G_P$. Note that this construction is in fact injective. The set of all these functions $g$ associated to $P\subseteq E_k$ that can be defined this way will be denote as $\mathfrak{G}_P$ in the following. In fact, the details of this construction won't really matter for the impatient reader since 
we will rapidly get rid of the $\mathfrak{G}_P$'s soon enough. Simply remember that $g$ is a kind of a "global permutation" made (uniquely) from all little permutations given on each face of a fixed graph $G_P, P\subseteq E_k$.
\\

Then the Relation \eqref{eq:bkhk}, together with what we just explained above, leads to:
\begin{equation}\label{eq:formulemoche}
-B(1\to k)=\sum_{g\colon [\![1,k]\!]\to\mathcal{T}([\![1,k]\!])}\,\,\sum_{\underset{g\in \mathfrak{G}_P}{P\subseteq E_k/}}(-1)^{\vert P\vert+1}\frac{1}{\prod_{e\in P\sqcup E_k^c}h(e)}\prod_{F \text{ face}}\left(C_{g_{|F}}(F)\!\!\prod_{e=(i,g_{|F}(i))}h(e)\right),
\end{equation}
where $g_{|F}$ stands for the restriction of $g$ on the face $F$ as explained in the previous paragraph, which hence is a plane permutation on $F$.

\paragraph{Simplifying terms via an involution.} To validate Theorem \ref{th:bk}, note that no $H_1$ terms should appear in the numerator: these must therefore cancel.
\\

Let us define an involution $\theta$ acting on the set of pairs $(g,P)$ such that $P\subseteq E_k$ and $g\colon [\![1,k]\!]\to \mathcal{T}([\![1,k]\!])\in\mathfrak{G}_P$.
 Assume further that there exists a pair $(i,j)$, with $j$ appearing in $g(i)$, belonging to $E_k$. We also require $i$ to be such that, if the edge $x_i\sim y_j$ lies in $P$, there exists another index $j'$ (possibly equal to $j$, should $j$ appear redundantly in $g(i)$) with $x_i\sim y_{j'}$ also in $P$. Set an arbitrary total order on $E_k$, and let $e$ be the minimal edge for this order, determined by the pair $(g,P)$ as described above.
\\

The map $\theta$ then simply adds $e$ to $P$ if it is absent, or removes it if present. The fact that $e$ lies in $E_k$ guarantees that $\theta$ is well-defined, and the construction clearly makes it an involution. Note also that our choice of $e$ ensures that adding or removing it preserves the plane bijectivity of the restriction of $g$ to each face of $G_{P\cup\{e\}}$ or $G_{P\backslash\{e\}}$. Let us now define
$$
w(P):=(-1)^{\vert P\vert+1}\frac{1}{\prod_{e\in P\sqcup E_k^c}h(e)}\prod_{F \text{ face}}\left(C_{g_{|F}}(F)\!\!\prod_{e=(i,g_{|F}(i))}h(e)\right),
$$
the weight associated with such a subset $P$ satisfying the above conditions relative to $g$. It remains to compute $w(\theta(P))$.
\\

First observe that, for fixed $g$, the products over $h$ do not change between $P$ and $\theta(P)$: adding $e$ to $P$ certainly introduces a new term in the denominator, but it cancels against a term arising in the numerator from the product over faces. It thus remains to examine the product of the $C$ terms, for which we require the former Lemma \ref{lem:csigma_face12}.\\

Thanks to this lemma — and above all to Equation \eqref{eq:csigma_face12} — we see that the product of $C$ terms in the definition of $w(P)$ is unaffected by adding or removing an edge from $P$. The only change is in the sign preceding it, since the number of edges shifts by $1$ under $\theta$. We have thus shown that
$$
w(\theta(P))=-w(P),
$$
and so a great many terms cancel in the sum \eqref{eq:formulemoche}! What remains are only the terms corresponding to functions $g$ whose image of each element contains at most $1$ or $2$ non-zero terms, always containing either $i$ or $i-1$ (or both) in $g(i)$. When $g(i)$ has $2$ non-vanishing elements, note that $P$ must contain an edge leaving $i$, or else $g$ would fail to belong to $\mathfrak{G}_P$.

\paragraph{Connection with plane permutations.}
As we have just seen, only a handful of functions $g$ survive in our sum, and the associated subsets $P\subseteq E_k$ then form the graphs of plane permutations, since $P$ cannot contain more than one pair of the form $(i,j)$ with $j\neq i,i-1$ appearing in $g(i)$. We may thus associate each such pair with a plane permutation, and rewrite Equation \eqref{eq:formulemoche} as:
\begin{equation}\label{eq:dsigma}
\boxed{-B(1\to k)=\sum_{\sigma\in\mathfrak{P}_k}D_{\sigma}(1\to k)\cdot\prod_{i=1}^k H_1(x_i;y_{\sigma}(i))^{-1},}
\end{equation}
where the $D_{\sigma}$'s remain to be determined. The goal of the next, and last, subsection is to show that $D_{\sigma}(1\to k)=C_{\sigma}(1\to k)$ for every $\sigma\in\mathfrak{P}_k$, which will establish Theorem \ref{th:bk}.

\subsection{Conclusion of the proof of Theorem \ref{th:bk}}

Now that we expressed $B$ as a sum over plane permutations exactly as $H_k$ is, it is time to go back to equation \eqref{eq:bh_bhmoinsun} in order to conclude.

\paragraph{Back to Equation \eqref{eq:bh_bhmoinsun}.} Going back and rewriting Equation \eqref{eq:bh_bhmoinsun} using our new formula \eqref{eq:dsigma} yield:
\begin{equation}\label{eq:Dsigma2}
\sum_{\sigma\in\mathfrak{P}_k}\mathcal{D}_{\sigma}(1\to k)\prod_{i=1}^kH_1(x_i;y_{\sigma}(i))=\sum_{\sigma\in\mathfrak{P}_k}\mathcal{D}_{\sigma}(1\to k)\frac{\prod_{i=1}^k H_1(x_i;y_i)\,H_1(x_i;y_{i-1})}{\prod_{i=1}^kH_1(x_i;y_{\sigma}(i))},
\end{equation}
where $\mathcal{D}_{\sigma}(1\to k):=C_{\sigma}(1\to k)-D_{\sigma}(1\to k)$.
A theorem proved in \cite{EO2005}, and reproved via the slice method in \cite{AB2025} and again in \cite{BEL26}, allows us to conclude. Let us recall it:
\begin{thm}[\cite{EO2005, AB2025, BEL26}]\label{th:h1param}
There exists a rational parametrization of $H_1$ of the form:
$$
H_1(x(z),y(z'))=\frac{E(x(z),y(z'))}{(x(z)-x(z'))(y(z)-y(z'))},
$$
where $x(z)$ and $y(z)$ are Laurent polynomials and $E$ is their resultant, hence itself a polynomial. Note also that $H_1(x(z),y(z'))$ cannot vanish since each zero of the resultant is compensated by the denominator.
\end{thm}
The careful reader will notice that the Laurent polynomials $x(z)$ and $y(z)$ are in fact the same as those appearing in Formula \eqref{eq:ak}, characterizing the alternating hypermaps $A_p, p\geqslant1$.\\

Now, parametrizing $x_i$ as $x(z_i)$ and $y_j$ as $y(z_j')$, we observe that $H_1(x(z_i);y(z_j'))$ has a pole at $z_i=z_j'$. Returning to the previous Equation \eqref{eq:Dsigma2} with this in mind, set $\sigma\in\mathfrak{P}_k$.
We see that, on the left-hand side, the term $\prod_{i=1}^kH_1(x(z_i);y(z'_{\sigma}(i)))$ is the only one with (simple) poles at $z_i=z'_{\sigma(i)},\forall i\in[\![1,k]\!]$. Indeed, the definitions of $C_{\sigma}(1\to k)$ and of $D_{\sigma}(1\to k)$ 
tell us that their only poles can exist only if two $x_i$'s (or two $y_i$'s) are equal, but never if $x_i$ equals one of the $y_j$'s. 
\\
Now look at the right-hand side and assume that $\sigma\neq \mathrm{Id},\gamma_k$. By the last remark of Theorem \ref{th:h1param}, the poles from the right-hand side are all the ones from the product $$\prod_{i=1}^kH_1(x(z_i);y(z'_i))\,H_1(x(z_i);y(z'_{i-1})),$$ and they are localised at $z_i=z'_i$ and $z_i=z_{i-1}$. But since $\sigma$ is neither $\mathrm{Id}$ nor $\gamma_k$, then there exists $i_0\in[\![1,k]\!]$ such that 
$\sigma(i_0)\notin\{i_0,i_0-1\}$. Finally, the term $\prod_{i=1}^kH_1(x(z_i);y(z'_{\sigma}(i)))$ is the only one on the left-hand side having poles localised at $z_i=z'_{\sigma(i)},\forall i$ and none of the terms of the right-hand side exhibit a pole at $z_{i_0}=z'_{\sigma(i_0)}$. Hence the coefficient in front of $\prod_{i=1}^kH_1(x(z_i);y(z'_{\sigma}(i)))$ must vanish, so 
$\mathcal{D}_{\sigma}(1\to k)=0$ or, equivalently $D_{\sigma}(1\to k)=C_{\sigma}(1\to k)$ for every $\sigma\neq \mathrm{Id}, \gamma_k$. 
\\

Now we only have to prove that $C_{\mathrm{Id}}(1\to k)=D_{\mathrm{Id}}(1\to k)$ and $C_{\gamma_k}(1\to k)=D_{\gamma_k}(1\to k)$. To do so, let us compare expand $H_k$ and $B_k$ as a series in the vertex marking variable $t$, which is equivalent to enumerate hypermaps with respectively $0$ and $1$ vertices. This gives, using Equation \eqref{eq:Hkdef}:
$$
H_k(1\to k;t)=\delta_{1,k}\,t^0+\frac{1}{\prod_{i=1}^k x_iy_i}\,t^1+O(t^2)
$$
and 
$$
B_k(1\to k;t)=0\,t^0+\frac{1}{\prod_{i=1}^k x_iy_i}\,t^1+O(t^2)
$$
since there exists only one hypermap with one vertex. Note that, for $k\geqslant2$, the $t^0$ and the $t^1$ terms from $H_k$ and $B_k$ match. Hence, extracting also these first two terms in Relations \eqref{eq:hkh1} and \eqref{eq:dsigma}, we obtain the relations:
$$
\forall k\geqslant2,\sum_{\sigma\in\mathfrak{P}_k}D_{\sigma}(1\to k)=[t^0]B_k(1\to k;t)=[t^0]H_k(1\to k;t)=\sum_{\sigma\in\mathfrak{P}_k}C_{\sigma}(1\to k),
$$
as well as
$$
\sum_{\sigma\in\mathfrak{P}_k}D_{\sigma}(1\to k)\sum_{i=1}^k\frac{1}{x_iy_{\sigma(i)}}=[t^1]B_k(1\to k;t)=[t^1]H_k(1\to k; t)=\sum_{\sigma\in\mathfrak{P}_k}C_{\sigma}(1\to k)\sum_{i=1}^k\frac{1}{x_iy_{\sigma(i)}}.
$$
Given the equalities $D_{\sigma}=C_{\sigma}$ for all $\sigma\neq \mathrm{Id},\gamma_k$ already established, only two equations remain, determining $D_{\mathrm{Id}}$ and $D_{\gamma_k}$ in terms of $C_{\mathrm{Id}}$ and $C_{\gamma_k}$. This system is clearly invertible, so the only possible solution is $D_{\mathrm{Id}}=C_{\mathrm{Id}}$ and $D_{\gamma_k}=C_{\gamma_k}$ — completing the proof of Theorem \ref{th:bk}!

\section{Some concluding remarks}

The purpose of this article was to describe hypermaps subject to general boundary condition, using a new tool we aim to develop further: accessibility. The idea is to mark a vertex of a hypermap and then consider its accessible component — that is, the set of vertices that can reach it. This approach lets us account for distances between points through the connection with slices, which is exactly what will be needed to go further and, for instance, build models of decorated random surfaces. \\
The only difficulty in our proofs lies in integrating the equations we derive, a step lacking any combinatorial interpretation, which makes the end of the argument more technical and less visual.
\\

We refer the reader interested in blobs — hypermaps in which every vertex is accessible — 
to the last section devoted to them. A proof similar to the one carried
 out in the second section of this article can also be given in that setting, 
 with a few additional subtleties concerning the orientation of the edges 
 considered in our decomposition into an accessible component. We chose not 
 to elaborate on this here so as not to overload the article but we will be very glad 
 to discuss it for further work purpose. Furthermore, a similar line of reasoning to that of the third section appears applicable to simple hypermaps — those in which no boundary vertex is incident to two corners lying outside the map — by drawing an analogy between such vertices and bridges. This would make for a fine continuation of this work.
\\

One remark that can nonetheless be made is that the coefficients
$C_{\sigma}$ appearing in the main equation of this corollary do not seem to
be combinatorial in themselves, even though the objects that lead to them
are. This is why we shall need some computations as well as an induction in
order to prove it, once Theorem \ref{th:dHk/dt} has been established, so as
to integrate the equation obtained. However, the accessibility method is not
meant to be difficult to integrate in every situation. We thus refer the
reader to the article \cite{ABL26}, which proposes to use this method to
treat the case of a fully alternating boundary, that is, alternating hypermaps where 
each boundary interval has length exactly one. We found there that the algebraicity of the
generating series is found more easily than here.
\\

Finally, the purpose of this method is not only to recover, by another means,
formulas that are already understood, but rather to offer a new approach to
the combinatorics of hypermaps. In particular, the connection with the slices
introduced in \cite{AB2025} and \cite{BOUTTIER_2014} makes it possible to compute
various asymptotics and, we hope, to generalize our approach to other models
such as the Potts model, the matrix chain model, the $O(n)$ loop model, the
Ising model, and so on. Since the bicolor case is the first step before
generalizing to a larger number of colors, this is why we are so keen to
strengthen the foundations of this first step before climbing the next ones.

\printbibliography
\end{document}